\documentclass[10pt]{amsart}

\usepackage[usenames,dvipsnames]{xcolor}
\usepackage[bookmarks,colorlinks=true,citecolor=OliveGreen,linkcolor=RoyalBlue]{hyperref}
\usepackage{amssymb,amsthm}
\usepackage[inline]{enumitem}
\usepackage{mathrsfs}
\usepackage{mathtools}
\usepackage{mathabx}
\usepackage{thmtools}
\usepackage[capitalise]{cleveref}		
\usepackage{bm}

\usepackage{lmodern}

\makeatletter
\def\thmt@refnamewithcomma #1#2#3,#4,#5\@nil{%
  \@xa\def\csname\thmt@envname #1utorefname\endcsname{#3}%
  \ifcsname #2refname\endcsname
    \csname #2refname\expandafter\endcsname\expandafter{\thmt@envname}{#3}{#4}%
  \fi
}
\makeatother

\declaretheorem[numberwithin=section]{theorem}
\declaretheorem[sibling=theorem]{proposition}
\declaretheorem[sibling=theorem]{corollary}

\declaretheorem[sibling=theorem,style=definition]{definition}

\declaretheorem[sibling=theorem,style=definition]{lemma}
\declaretheorem[sibling=theorem,style=remark]{remark}

\declaretheorem[sibling=theorem,style=remark]{notation}

\newcounter{claimCounter}[theorem]

\newcounter{subClaimCounter}[claimCounter]

\newlist{equivalent}{enumerate}{1}
\setlist[equivalent,1]{label=\textup{(\arabic*)}}

\newlist{sublemma}{enumerate}{1}
\setlist[sublemma,1]{label=\textup{(\alph*)}}

\newlist{sublemma*}{enumerate*}{1}
\setlist[sublemma*,1]{label=\textup{(\alph*)},afterlabel=\hspace{5pt}}

\newlist{orderedlist}{enumerate}{1}
\setlist[orderedlist,1]{label=\textup{(\roman*)}}

\newlist{orderedlist*}{enumerate*}{1}
\setlist[orderedlist*,1]{label=\textup{(\roman*)},afterlabel=\hspace{3pt}}

\newcommand{\seq}[1]{{\left\langle{#1}\right\rangle}}
 
\newcommand{\rest}[1]{\! \upharpoonright\!{#1}} 

\newcommand{\conc}{\hat{\,\,}}
\newcommand{\andd}{\,\,\,\&\,\,\,}

\newcommand{\Iff}{\,\,\Longleftrightarrow\,\,}

\newcommand{\smallseq}[1]{{\langle{#1}\rangle}}

\DeclareMathOperator{\dom}{dom}

\DeclareMathOperator{\cf}{cf}

\newcommand{\w}{\omega}
\newcommand{\s}{\sigma}
\newcommand{\vphi}{\varphi}
\renewcommand{\epsilon}{\varepsilon}

\renewcommand{\le}{\leqslant}
\renewcommand{\ge}{\geqslant}

\renewcommand{\preceq}{\preccurlyeq}
\renewcommand{\succeq}{\succcurlyeq}

\newcommand{\Int}{\mathbb Z}

\newcommand{\PP}{\mathbb P}
\newcommand{\RR}{\mathbb R}

\newcommand{\task}{\mathfrak{t}}
\newcommand{\HH}{\mathcal{H}}
\newcommand{\FF}{\mathscr{F}}
\newcommand{\PPP}{\mathcal{P}}
\newcommand{\starr}{\!*\!}
\newcommand{\QQ}{\mathbb Q}
\renewcommand{\SS}{\mathbb S}

\DeclareMathOperator{\otp}{otp}

\DeclareMathOperator{\NS}{NS}
\newcommand \compembeded{\lessdot}

\newcommand{\GCH}{\textup{GCH}}

\newcommand{\force}{\Vdash}

\newcommand{\lth}[1]{\left|#1\right|}
\newcommand{\smalllth}[1]{|#1|}

\newcommand{\cohere}{\lessdot}

\newcommand{\DiagUnion}{\nabla}

\newcommand{\BB}{\mathscr{B}}
\newcommand{\II}{\mathbf{I}}
\newcommand{\ff}{\mathbf{f}}

\DeclareMathOperator{\Cof}{Cof}

\DeclareMathOperator{\supp}{supp}

\newcommand{\uu}{\mathbf{u}}

\newcommand{\lllambda}{<\!\!\lambda}
\newcommand{\UU}{\mathscr{U}}

\DeclareMathOperator{\cb}{cb}

\newcommand{\TaskAxiom}[1]{\textup{\textrm{TaskAx}}_{#1}}

\newcommand{\UP}{\textup{\scriptsize UP}}

\DeclareMathOperator{\Ext}{Ext}

\title{Many forcing axioms for all regular uncountable cardinals}

\author[N.~Greenberg]{Noam Greenberg} 
\address[N.~Greenberg]{School of Mathematics, Statistics and Operations Research, Victoria University of Wellington, Wellington, New Zealand}
\email{greenberg@msor.vuw.ac.nz}
\urladdr{\url{http://homepages.mcs.vuw.ac.nz/~greenberg/}}

\author[S.\ Shelah]{Saharon Shelah}
\address{The Hebrew University of Jerusalem \\ Einstein Institute of
Mathematics \\ Edmond J. Safra Campus, Givat Ram \\ Jerusalem 91904,
Israel}
\address{Department of Mathematics \\ Hill Center-Busch Campus \\ Rutgers,
The State University of New Jersey \\ 110 Frelinghuysen Road
\\ Piscataway, NJ 08854-8019 USA}
\email{shelah@math.huji.ac.il}

\thanks{Greenberg was partially supported by a Rutherford Discovery Fellowship and a Lady Davis visiting professorship. Shelah's research partially supported by the German-Israeli Foundation for scientific research and development, grant no.: I-706054.6/2001; and by the Israel Science Foundation, grant no.: ISF 1838/19. Paper 832 on Shelah's list.}

\begin{document}

\begin{abstract}
		A central theme in set theory is to find universes with extreme, well-understood behaviour. The case we are interested in is assuming GCH and having a strong forcing axiom of higher order than usual. Instead of ``every suitable forcing notion of size~$\lambda$ has a sufficiently generic filter''	we shall say ``for every suitable method of producing notions of forcing based on a given stationary set, there is such a suitable stationary set~$S$ and sufficiently generic filters for the notion of forcing attached to~$S$''. Such notions of forcing are important for Abelian group theory, but this application is delayed for a sequel. 
\end{abstract}

\maketitle

\section{introduction} 

One of the original motivations for the work presented in this paper is to show the consistency of the failure of singular compactness for properties such as being a Whitehead group ($\Ext(G,\Int)=0$). Under $V=L$, for example, singular compactness for being a Whitehead group holds, but this is because $V=L$ implies that Whitehead groups are all free \cite{Sh:44}, and singular compactness for free groups is a ZFC theorem \cite{Sh:SingularCompactness}. The question whether $V=L$ is necessary arises from work in \cite{Eklof:80}; see also \cite{EkFuSh:352,Hodges,Sh:266}. A similar question about singular compactness of a property $\Ext(G,T)=0$ where~$T$ is a torsion group was asked by Str\"{u}ngmann, following \cite[Prop.2.6]{Str:02}. 

These problems are closely related to so-called uniformisation principles (see~\cite{EklofMekler:book}), and in many cases are even equivalent to them. The first work along these lines is~\cite{Sh:64}, where it is proved that GCH is consistent with diamond holding at some stationary subset of~$\w_1$, while failing at others (indeed on the latter, some uniformisation principle, usually derived from MA, holds). For more see \cite{Sh:587,Sh:667}. 

In this paper we present an axiom which guarantees, among other consequences, instances of uniformisation. What we term the ``task axiom'' for a regular cardinal~$\lambda$ ensures that mutually competing principles, such as diamonds and uniformisation, all hold at (necessarily disjoint) stationary sets. 

A feature of this axiom is that it entails~$\Pi_4$, rather than~$\Pi_2$ statements: to satisfy a ``task'', we need not only a sufficiently generic filter for an appropriate notion of forcing, but also many filters for subsequent notions of forcing which are determined by the first filter. The two examples we present in this paper are exact diamonds (\cref{prop:task_axiom_and_diamond}) and uniformisation (\cref{prop:uniformisation}). For example, for the former, the first notion of forcing adds a stationary set and a diamond sequence on that set; subsequent notions of forcing ensure that this diamond sequence is exact, by adding functions which are not guessed by less than the full diamond sequence. 

Much of the technical work is involved in showing that the stationary set, added by the primary notion of forcing of a task (in the example above, the set on which the diamond sequence is supported), remains stationary after adding the subsequent generic filters. We also need to ensure that different tasks do not interfere with each other; that is more easily achieved, by ensuring that distinct tasks are met on disjoint stationary sets. Overall, this requires a ``niceness condition'' which is part of the definition of the task axiom. The definition of this condition uses a notion of closure of a forcing partial ordering on a given fat set, related to the notions of $S$-completeness defined in \cite{Sh:64,Sh:587,Sh:667}. 

In a planned sequel, the second author intends to use the task axiom to show the consistency of the failure of singular compactness for the classes of groups mentioned above.

\subsection{The contents of the paper} 

Tasks (\cref{def:task}) are defined in \cref{sec:tasks_and_the_task_iteration}; the task axiom (\cref{def:lambda_task_axiom}) is stated in \cref{sec:the_task_axiom}. Beforehand, we develop the tools required to formulate tasks and to work with them. Throughout, we fix a regular uncountable cardinal~$\lambda$, and assume that $\GCH$ holds below~$\lambda$. 

To start, we define two notions of completeness:
\begin{itemize}
	\item explicitly $S$-closed notions of forcing (\cref{def:explicitly-S-closed}); 
	\item $S$-sparse names for subsets of~$\lambda$ (\cref{def:sparse_choice_sparse_closure}). 
\end{itemize}
The former attaches ordinals to conditions, and states that increasing sequences of conditions whose associated oridnals converge to points in~$S$ have upper bounds. The latter is similar, except that we now associate to conditions initial segments of a subset of~$\lambda$ in the generic extension, and state that increasing sequences of conditions which determine a closed set disjoint from the subset in question have upper bounds. 

We then turn to develop machinery that will help us work with forcing iterations. The general situation is a $\lllambda$-support iteration $\seq{\PP_\zeta,\QQ_\zeta}_{\zeta < \xi}$, one of which we will use to show consistency of the task axiom. In the intended application, some of the notions of forcing~$\QQ_\zeta$ add new stationary sets~$W_\zeta$, on which we will want to fulfill some ``task''; that is done by some of the subsequent notions of forcing~$\QQ_{\eta}$ appearing (cofinally) later in the iteration. The typical task will state the existence of a stationary set~$S$ and possibly some associated object (for example a diamond sequence), for which some $\Pi_2(\HH_\lambda)$ fact holds (for example, the diamond sequence is exact).  

For simplicity of notation, we will set $W_{\eta} = \emptyset$ for the notions~$\QQ_{\eta}$ which do not start a new task. Thus, for many of our definitions and lemmas, we do not need to differentiate between the two kinds of $\QQ_{\zeta}$'s. We will ensure that the sets~$W_{\zeta}$ are pairwise disjoint (modulo clubs). 

In \cref{sec:closure_and_sparseness} we show (\cref{cor:iteration_up_to_lambda_plus:existence_of_associated_sequence}) how to construct, for each $\zeta \le \xi$, a sparse $\PP_\zeta$-name for the least upper bound of the sets $W_{\zeta'}$ (for $\zeta'<\zeta$), again modulo clubs. These names, and how they cohere with each other, will be one of the main tools we use in the analysis of tasks, ensuring that they can be adequately fulfilled. 

In \cref{sec:tasks_and_the_task_iteration} we define tasks. To motivate the definitions, before we define the task axiom, we give a prototypical example of an iteration, of length~$\lambda^+$, attempting to fulfill tasks. We formulate the example of exact diamonds and show that these exist in the generic extension (\cref{prop:first_iteration:diamonds_in_extension}). 

To formulate the task axiom, we need to bear in mind that some tasks may be too ambitious; an attempt to fulfill them, along with other tasks, will result in failure, for example, in a hoped-for stationary set not being stationary anymore. Thus we need to define ``correctness conditions'' for tasks (\cref{def:correctness_condition}). It is only reasonable to expect that a task be fulfilled if a typical forcing iteration preserves this correctness condition; this is the notion of niceness for tasks (\cref{def:nice_task}), with which we can then state the task axiom. 

Unfortunately, the iteration of length~$\lambda^+$ given in \cref{sec:tasks_and_the_task_iteration} may be too short to expose those tasks which are not nice, i.e., to witness the failure of the correctness condition after some iteration. The consistency of the task axiom is then proved by a longer iteration (still of length $<\lambda^{++}$), with sufficient closure. At steps of cofinality~$\lambda^+$, we need auxiliary notions of forcing to ensure that the sequence of sets~$W_\zeta$ have upper bounds; we develop this machinery (\cref{prop:obtaining_tau_for_long_iteration}) in the beginning of \cref{sec:the_task_axiom}. 

Finally, as a second example, we show that the task axiom implies instances of uniformisation (\cref{prop:uniformisation}). 

We remark that the construction can be modified to ensure the preservation of suitable large cardinals; we do not pursue this topic in this paper.

\subsection{Notation and terminology} 
\label{sub:notation_and_terminology}

We list some of the notation that we use. We follow the Israeli convention for extension in notions of forcing; $p\le q$ means that~$q$ extends~$p$. Complete embeddings $\PP\cohere \QQ$ of forcing notions will always be accompanied with a complete projection $p\mapsto p\rest{\PP}$ from~$\QQ$ to~$\PP$. For typographical convenience, for a notion of forcing~$\PP$ and a definable set or class~$X$, we let $X(\PP)$ denote the interpretation of~$X$ in the Boolean valued model $V^\PP$. So for example $V(\PP) = V^\PP$, $\PPP(\lambda)(\PP) = (\PPP(\lambda))^{V^\PP}$ is the collection of names $a\in V^\PP$ such that $\force_\PP a\subseteq \lambda$, etc. When we say that a statement~$\vphi$ holds in $V(\PP)$, we mean that every condition forces~$\vphi$. 

For a binary string $\s$ (a function from an ordinal into 2), we let $\lth{\s} = \dom \s$ denote the length of~$\s$. We will at times be imprecise and identify sets with characteristic functions; that is, we identify $\s\colon \alpha \to 2$ with $\left\{ \beta < \alpha \,:\,  \s(\beta)=1 \right\}$, when there is no danger that we forget $\alpha = \lth{\s}$. we write $\s\preceq \tau$ to indicate that~$\tau$ extends~$\s$, that is, $\s = \tau\rest{\lth{\s}}$.

$\Cof(\theta)$ denotes the class of ordinals of cofinality $\theta$; similarly we use $\Cof(\le \theta)$ etc. For a cardinal~$\lambda$, we let $\Cof_\lambda(\theta)$ denote $\lambda\cap \Cof(\theta)$. 

We let $\+H_\chi$ denote the collection of sets whose transitive closure has size $<\chi$. 

\subsection{The underlying hypothesis} 

Throughout this paper, $\lambda$ denotes a regular uncountable cardinal, and we assume that $\GCH$ holds below~$\lambda$.

\subsection{Approachable ordinals} 
\label{sub:approachable_ordinals}

In this section we isolate a tool that will allow us to make use of closure conditions on fat sets. The technique we use was introduced in \cite{Sh:108}, where it is used to show that there are many ``approachable'' ordinals below~$\lambda$ (so a large subset of a given fat subset of~$\lambda$ will be in $\check I[\lambda]$). Recall that a subset~$S$ of~$\lambda$ is \emph{fat} iff for every regular $\theta <\lambda$, for every club~$C$ of~$\lambda$, $S\cap C$ contains a closed subset of order-type $\theta+1$. If~$S$ is fat then for all regular $\theta<\lambda$, $S\cap \Cof(\theta)$ is stationary in~$\lambda$. 

For the following lemma, recall that a \emph{$\lambda$-filtration of the universe} is an increasing and continuous sequence $\bar N = \seq{N_{\gamma}}_{\gamma<\lambda}$ such that for some large~$\chi$, for all $\gamma<\lambda$,
\begin{orderedlist}
	\item $N_{\gamma}$ is an elementary submodel of $\+H_\chi$;
	\item $|N_{\gamma}|<\lambda$;
	\item $\gamma\subseteq N_{\gamma}$; and
	\item $\bar N\rest{(\gamma+1)}\in N_{\gamma+1}$. 
\end{orderedlist}

\begin{lemma} \label{lem:approachable_ordinals}
	Suppose that $S\subseteq \lambda$ is fat. For every regular cardinal $\theta<\lambda$ and every set~$X$ with $|X|<\lambda$ there is a $\lambda$-filtration~$\bar N$ such that $X\subseteq N_{0}$, and a set $D\subset \lambda$ such that:
	\begin{orderedlist}
		\item $\otp(D)= \theta+1$;
		\item The closure $\overline D$ of~$D$ is a subset of~$S$;
		\item For all $\gamma\in D$, $N_{\gamma}\cap \lambda = \gamma$;
		\item For all $\gamma\in D$ other than $\max D$, $D\cap \gamma\in N_{\gamma+1}$. 
	\end{orderedlist}
\end{lemma}

Note that the set~$X$ needn't be a subset of~$\lambda$; we will have $X\in \HH_\chi$ for some large~$\chi$, of which the~$N_\gamma$ will be elementary submodels.

\begin{proof}
	There are two cases. The easier one is when~$\lambda$ is not the successor of a singular cardinal, that is, it is either the successor of a regular cardinal, or inaccessible. In that case we can build the filtration $\seq{N_{\gamma}}$ inductively, requiring that for all $\gamma<\lambda$, $[\gamma]^{<\theta}\subset N_{\gamma+1}$. This we can do because $|\gamma|^{<\theta}<\lambda$, as in either case there is a regular cardinal $\kappa<\lambda$ such that $|\gamma|,\theta\le \kappa$, and for all regular $\kappa<\lambda$, $\kappa^{<\kappa}=\kappa$. We then let $E = \left\{ \gamma<\lambda \,:\,  N_{\gamma}\cap \lambda = \gamma  \right\}$; this is a club of~$\lambda$, so $S\cap E$ contains a closed subset~$D$ of order-type $\theta+1$. If $\gamma\in D$ and $\gamma\ne \max D$ then $|D\cap \gamma|<\theta$, so is an element of $N_{\gamma+1}$. 

	\smallskip
	
	We now suppose that~$\lambda = \mu^+$ where~$\mu$ is singular; let $\kappa=\cf(\mu)$. In this case we may have $|\gamma|^{<\theta}= \lambda$, so we need a finer approach. We know that $\theta<\mu$. 

	Let $\seq{\mu_\xi}_{\xi<\kappa}$ be a sequence of cardinals increasing to~$\mu$. For all $\alpha\in [\mu,\lambda)$ fix a partition $\{A^\alpha_\xi\,:\, \xi<\kappa\}$ of~$\alpha$ such that $|A^\alpha_\xi| = \mu_\xi$. We build $\seq{N_\gamma}_{\gamma<\lambda}$ such that $|N_\gamma|= \mu$ and ensure that all bounded subsets of~$\mu$ are in~$N_0$. We also put the map $(\alpha,\xi)\mapsto A^\alpha_\xi$ into~$N_0$. We define~$E$ as above. 

	We may assume that $\theta>\kappa$. By $\GCH$ below~$\lambda$, we know that $(2^\theta)^+<\lambda$, and so we can find a closed set $D^*\subset S\cap E$ of size $(2^\theta)^+$. By the Erd\"{o}s-Rado theorem we can find $D\subset D^*$ of order-type $\theta+1$ (in fact of size $\theta^+$) such that for all $\gamma<\delta$ from~$D$, $\gamma\in A^\delta_{\xi^*}$ for some fixed $\xi^*<\kappa$. The set~$D$ may not be closed but its closure is a subset of~$D^*$, and so of $S\cap E$. For all $\gamma\in D$, $D\cap \gamma\subseteq A^\gamma_{\xi^*}$. Since $A^\gamma_{\xi^*}\in N_{\gamma+1}$, and this set is bijective with $\mu_{\xi^*}$, every subset of $A^\gamma_{\xi^*}$ is in $N_{\gamma+1}$. So $D\cap \gamma\in N_{\gamma+1}$. 
\end{proof}


\section{Closed notions of forcing and sparse names} 
\label{sec:closure_and_sparseness}


\subsection{Explicit $S$-closure} 
\label{sub:explicit_S_closure}

\begin{definition} \label{def:explicitly-S-closed}
	Let $S\subseteq \lambda$. A notion of forcing~$\PP$ is \emph{explicitly $S$-closed} if there is a function $\delta \colon \PP\to \lambda$ satisfying:
	\begin{orderedlist}
		\item $p\le q$ implies $\delta(p)\le \delta(q)$;
		\item For every $\gamma<\lambda$, the collection of conditions~$p$ with $\delta(p)\ge \gamma$ is dense in~$\PP$;
		\item whenever $\bar p = \seq{p_i}_{i<i^*}$ is an increasing sequence of conditions from~$\PP$, with $i<j<i^{*}$ implying $\delta(p_{i})< \delta(p_{j})$, and $\alpha = \sup_{i<i^*} \delta(p_i) \in S$, then $\bar p$ has an upper bound~$p^{*}$ in~$\PP$ with $\delta(p^{*}) = \alpha$. 
	\end{orderedlist}
\end{definition}

We call such an upper bound~$p^*$ an \emph{exact upper bound} of~$\bar p$.

\begin{lemma} \label{lem:explicit_closure_and_club_equivalence}
	If $S$ and~$S'$ are equivalent modulo the club filter, then a notion of forcing~$\PP$ is explicitly $S$-closed if and only if it is explicitly $S'$-closed. \qed	
\end{lemma}

\begin{proof}
	If $\PP$ is explicitly $S$-closed, as witnessed by~$\delta$, and~$C$ is a club, then by replacing~$\delta(p)$ by $\sup (C\cap \delta(p))$, we may assume that~$\delta$ takes values in~$C$. 
\end{proof}

\begin{proposition} \label{prop:adding_shorter_sequences_with_fatness_and_explicit_closure}
	Suppose that
		 $S\subseteq \lambda$ is fat, and that 
		$\PP$ is explicitly $S$-closed.
	Then~$\PP$ is $\lllambda$-distributive, and~$S$ is fat in $V(\PP)$. 
\end{proposition}

\begin{proof}
	Let $\theta<\lambda$ be a regular cardinal, and let $\left\{ U_i \,:\,  i<\theta \right\}$ be a family of dense open subsets of~$\PP$; let $p_0\in \PP$. Let $\bar N = \seq{N_\gamma}_{\gamma<\lambda}$ and~$D$ be given by \cref{lem:approachable_ordinals}, with $\PP,p_0, \seq{U_i}\in N_0$.

	We enumerate~$D$ as $\seq{\gamma_i}_{i\le \theta}$ and define an increasing sequence $\bar p = \seq{p_i}_{i\le \theta}$ (starting from $p_0$) of conditions from~$\PP$. We fix a well-ordering $\le^*_\PP$ of~$\PP$ which is an element of~$N_0$, and define the sequence of conditions as follows:
	\begin{sublemma}
		\item Given $p_i$ for $i<\theta$, $p_{i+1}$ is the $\le^*_\PP$-least $p\in \PP$ extending~$p_i$ such that $\delta(p)\ge \gamma_i$ and $p\in U_i$. 
		\item For limit $i\le \theta$, $p_i$ is the $\le^*_\PP$-least upper bound of $\bar p\rest{i}$. 
	\end{sublemma}

	Of course we need to argue that this construction is possible, which means showing that at a limit step $i\le \theta$, $\bar p \rest{i}$ has an upper bound. First we observe that if the construction has been performed up to and including step $i+1$, then the sequence $\seq{p_j}_{j\le i+1}$ is definable from $\PP,p_0, \seq{U_i}, \le^*_\PP$ and the sequence $\seq{\gamma_j}_{j\le i}$. These parameters are all in $N_{\gamma_i+1}$, so $p_{i+1}\in N_{\gamma_i+1}$. It follows that $\delta(p_{i+1})< \gamma_{i+1}$. 

	Suppose now that $i\le \theta$ is a limit ordinal and that the sequence $\seq{p_j}_{j<i}$ has been successfully defined. Then as $\gamma_j \le \delta(p_{j+1})< \gamma_{j+1}$ for all $j<i$, we have $\sup_{j<i}\delta(p_j) = \sup_{j<i} \gamma_j$ is a limit point of~$D$, and hence is in~$S$; so the closure assumption says that $\bar p\rest{i}$ has an upper bound in~$\PP$, whence~$p_i$ is defined and the construction can proceed. Then~$p_\theta$ is an extension of~$p_0$ in $\bigcap_i U_i$. 

	\smallskip
	
	To show that~$S$ is fat in $V(\PP)$, again let $\theta<\lambda$ be regular, and let $C\in V(\PP)$ be a club of~$\lambda$. We perform a similar construction, this time with $p_{i+1}$ forcing some $\beta_i>\gamma_i$ into~$C$. As $p_{i+1}\in N_{\gamma_i+1}$, we have $\beta_i < \gamma_{i+1}$. The final condition forces that the limit points of~$D$ are all in~$C$, and so form a closed set of order-type $\theta+1$ in $S\cap C$. 
\end{proof}

Note that if $S'\subseteq S$ then $\PP$ is also explicitly $S'$-closed, and so if~$S'$ is fat, then it is fat in $V(\PP)$. The argument for the case $\theta = \w$ shows that if $S'\subseteq S\cap \Cof(\aleph_0)$ is stationary, then it remains stationary in $V(\PP)$. 


\smallskip

We introduce terminology: \begin{definition}
	We say that~$\PP$ is \emph{explicitly closed outside~$S$} if it is explicitly $(\lambda\setminus S)$-closed.
\end{definition} 

\subsubsection{Strategic closure} 

We remark that we can generalise \cref{def:explicitly-S-closed}:

\begin{definition} \label{def:strategic-S-closure}
	Let $S\subseteq \lambda$. A notion of forcing~$\PP$ is \emph{strategically $S$-closed} if the ``completeness'' player has a winning strategy in the following game, of length $\le\lambda$: the players alternate choosing conditions forming an increasing sequence $\seq{p_i}$ from~$\PP$ (the incompleteness player chooses first; the completeness player chooses at limit steps), and also choosing ordinals below~$\lambda$ which form an increasing and continuous sequence $\seq{\epsilon_i}$ (so at limit stages there is no freedom in the choice of ordinal). 
	\begin{itemize}
		\item The incompleteness player loses at a limit stage~$i$ if $\epsilon_i \notin S$; 
		\item The completeness player loses at a limit stage~$i$ if $\epsilon_i \in S$ and $\seq{p_j}_{j<i}$ does not have an upper bound in~$\PP$. 
	\end{itemize}
	If the play lasts~$\lambda$ moves then the completeness player wins.
 \end{definition}

Note that if $S = \lambda$ then there is no need to choose ordinals; it is then the usual notion of $\lambda$-strategic closure. Also note that if~$S$ and~$S'$ are equivalent modulo the club filter, then~$\PP$ is strategically $S$-closed if and only if it is strategically $S'$-closed. 



\begin{lemma} 
	If $\PP$ is explicitly $S$-closed then it is strategically $S$-closed. \qed
\end{lemma}

\Cref{prop:adding_shorter_sequences_with_fatness_and_explicit_closure} holds with the weaker hypothesis that~$\PP$ is strategically $S$-closed. The proof is modified so that the sequence $\bar p$ is part of a play $(\bar p,\bar \epsilon)$ in which the completeness player follows her winning strategy; the incompleteness player chooses conditions in the desired dense sets, and plays the ordinal $\epsilon_i = \gamma_i$. The completeness player's response is inside $N_{\gamma_i+1}$, so the ordinal she plays is below~$\gamma_{i+1}$.

\begin{remark} \label{rmk:startegic_and_exmplicity}
		If~$\PP$ is strategically $S$-closed then there is some explicitly~$S$-closed notion of forcing~$\PP'$ and a complete projection from~$\PP'$ onto~$\PP$: fixing a strategy~$\mathfrak{s}$ witnessing strategic closure of~$\PP$, we let~$\PP'$ consists of all plays in which the completeness player follows~$\mathfrak{s}$ and which have a last move, made by the incompleteness player. This gives an alternative proof of \cref{prop:adding_shorter_sequences_with_fatness_and_explicit_closure} when~$\PP$ is strategically $S$-closed. 
\end{remark}

\subsection{Adding sparse sets} 
\label{sub:adding_sparse_sets}

\begin{definition} \label{def:explicitly_adding_a_set_with_initial_segments}
	Let~$\PP$ be a notion of forcing. An \emph{explicit $\PP$-name for a subset of~$\lambda$} is a partial map~$\s$ from $\PP$ to $2^{<\lambda}$ satisfying:
	\begin{orderedlist}
		\item If $p,q\in \dom \s$ and~$q$ extends~$p$ then $\s(q)\succeq \s(p)$;
		\item $\dom \s$ is dense in~$\PP$; in fact, for every $\gamma<\lambda$, the set of conditions $p\in \dom \s$ with $\lth{\s(p)}\ge \gamma$ is dense in~$\PP$.\footnote{Recall that $|\tau| = \dom \tau$ is the length of the string~$\tau$.}
	\end{orderedlist}
\end{definition}

That is, we insist that all bounded initial segments of the set named are determined in~$V$. To avoid confusion, we denote by $W_\s$ the actual $\PP$-name of the resulting subset of~$\lambda$: $(p,\alpha)\in W_\s$ if $\alpha<\lth{\s(p)}$ and $\s(p)(\alpha)=1$. In other words, if $G\subset \PP$ is generic, then 
\[
	 W_{\s}[G] = \left\{ \alpha<\lambda \,:\,  (\exists p\in G\cap \dom \s) \,\,\s(p)(\alpha)=1\right\}.
\]
Of course, if~$\PP$ is $\lllambda$-distributive then every subset of~$\lambda$ in $V(\PP)$ has an explicit name. We remark that it is not difficult to extend~$\s$ to be defined on all of~$\PP$, but we will naturally work with dense subsets of orderings in a way that makes this formulation more convenient.

\smallskip

Suppose that~$\s$ is an explicit $\PP$-name for a subset of~$\lambda$. For an increasing sequence $\bar p =\seq{p_i}_{i<i^*}$ of conditions from $\dom \s$ (where $i^*<\lambda$), we let 
\[
	\s(\bar p) = \bigcup_{i<i^*} \s(p_i). 
\]

An \emph{exact sparse upper bound} (with respect to~$\s$) for a sequence $\bar p$ is an upper bound~$p$ of~$\bar p$ in $\dom \s$ such that $\s(p) = \s(\bar p)\conc 0$. 

\begin{definition} \label{def:sparse_choice_sparse_closure}
	Let $\s$ be an explicit $\PP$-name for a subset of~$\lambda$.
	\begin{sublemma}
	\item	 A \emph{choice of sparse upper bounds} is a partial function~$f$, defined on a collection of increasing sequences of conditions from $\dom \s$, such that for each $\bar p\in \dom f$, $f(\bar p)$ is an exact sparse upper bound of~$\bar p$. 
	
	\item	Two increasing sequences $\bar p = \seq{p_i}_{i<i^*}$ and $\bar q = \seq{q_j}_{j<j^*}$ are \emph{co-final} if for all~$i$ there is a~$j$ such that $p_i\le q_j$ and vice-versa. We say that~$f$ is \emph{canonical} if whenever $\bar p$ and~$\bar q$ are co-final and $f(\bar p)$ is defined, then $f(\bar q)$ is defined as well, and $f(\bar q) = f(\bar p)$. 
	
	\item	We say that an increasing sequence $\bar p = \seq{p_i}_{i<i^*}$ is a \emph{sparse sequence} (for $\s$ and~$f$) if for all $i<j<i^*$, $\lth{\s(p_i)}< \lth{\s(p_j)}$, and for all limit $i<i^*$, $f(\bar p\rest{i})$ is defined and $p_i = f(\bar p\rest{i})$. 
	\end{sublemma}
\end{definition}

A sparse sequence $\bar p = \seq{p_i}_{i\le i^*}$ determines a closed set disjoint from $W_\s$, namely the set 
\[
	\left\{ \lth{\s(\bar p\rest{j})} \,:\,  j\le i^* \text{ limit} \right\};
\]
this is a set in~$V$, and $p_{i^*}$ forces that this set is disjoint from~$W_\s$ (which note is not in~$V$). 

\begin{definition} \label{def:sparse_name}
		Let $S\subseteq \lambda$. An \emph{$S$-sparse $\PP$-name} (for a subset of~$\lambda$) is a pair $(\s,f)$ consisting of an explicit $\PP$-name~$\s$ for a subset of~$\lambda$, and a choice~$f$ of sparse upper bounds for~$\s$, satisfying:
		\begin{itemize}
			\item For any sparse sequence~$\bar p$ for~$\s$ and~$f$ (of limit length), if $\lth{\s(\bar p)}\in S$, then $f(\bar p)$ is defined. 
		\end{itemize}
\end{definition}

\begin{notation}
	If $(\s,f)$ is a sparse name then we usually denote~$f$ by $\cb_\s$ (standing for ``canonical bound''). We use~$\s$ to also denote the pair $(\s,\cb_\s)$.
\end{notation}

\begin{remark} \label{rmk:countable_closure}
	Suppose that $\s$ is $S$-sparse. Any increasing $\w$-sequence~$\seq{p_n}$ of conditions from $\dom \s$ (with $\s(p_n)$ strictly increasing) is, vacuously, sparse for~$\s$, and so if $\lth{\s(\bar p)}\in S$ then $\cb_\s(\bar p)$ is defined. 
\end{remark}

\begin{lemma}  \label{lem:sparse_implies_strategic}
	If there is an $S$-sparse $\PP$-name then~$\PP$ is strategically $S$-closed.
\end{lemma}

\begin{proof}
 A winning strategy for the completeness player ensures that at limit steps, we use the canonical choice of a sparse upper bound, thus ensuring that the sequence of conditions played by the completeness player is sparse. The ordinal played at successor step~$i$ is $\lth{\s(p_i)}$. 
\end{proof}

\begin{lemma} \label{lem:sparseness_invariant_under_clubs}
	If~$\s$ is an $S$-sparse $\PP$-name, and $S'$ is equivalent to~$S$ modulo the club filter on~$\lambda$, then there is an $S'$-sparse $\PP$-name~$\s'$ such that in $V(\PP)$, $W_{\s} = W_{\s'}$.\footnote{Recall that this means that every condition forces that $W_\s = W_{\s'}$.}
\end{lemma}

\begin{proof}
	We extend the argument of \cref{lem:explicit_closure_and_club_equivalence}.	Suppose that $S\cap C = S'\cap C$ for a club~$C$. Determine that $\dom \s' = \dom \s$. For $p\in \dom \s$, let $\alpha = \sup (C\cap \lth{\s(p)})$. If $\alpha = \lth{\s(p)}$ then let $\s'(p) = \s(p)$. Otherwise let $\s'(p) = \s(p)\rest{(\alpha+1)}$. Suppose that $\bar p = \seq{p_i}_{i<i^*}$ is an increasing sequence of conditions with $i<j<i^*$ implying $\lth{\s'(p_i)}< \lth{\s'(p_j)}$. So $C\cap [\lth{\s'(p_i)}, \lth{\s'(p_j)})\ne \emptyset$. It follows that if $i^*$ is a limit then $\s(\bar p) = \s'(\bar p)$. We therefore use the same choice of canonical sparse upper bounds; if $\bar p$ is a sparse sequence for~$\s'$, then it is a sparse sequence for~$\s$. 
\end{proof}

\begin{proposition} \label{prop:fatness_preserved_by_explictly_adding_subset}
	Suppose that 
		$S\subseteq \lambda$ is fat, and 
		that~$\s$ is an $S$-sparse $\PP$-name.
	Then $\PP$ is $\lllambda$-distributive, and~$S\setminus W_{\s}$ is fat in~$V(\PP)$. 
\end{proposition}

And as above, the proof also shows that for all stationary $S'\subseteq S\cap \Cof(\aleph_0)$, $S'\setminus W_{\s}$ is stationary in~$V(\PP)$.

\begin{proof}
	By \cref{lem:sparse_implies_strategic}, we know that $\PP$ is $\lllambda$-distributive. 

	To see that $S\setminus W_{\s}$ is fat, we mimic the proof of \cref{prop:adding_shorter_sequences_with_fatness_and_explicit_closure}, ensuring that the sequence of conditions $\bar p$ is sparse, with $\lth{\s(p_{i+1})}> \gamma_i$. At limit steps we take the canonical sparse upper bound. Again $p_{i+1}\in N_{\gamma_i+1}$ forces some $\beta_i>\gamma_i$ into the club~$C$. The limit points of~$D$ are disjoint from~$W_{\s}$. 
\end{proof}


\begin{lemma} \label{lem:closed_implies_sparsely_closed}
	Suppose that~$\PP$ is explicitly $S$-closed. Then there is an $S$-sparse $\PP$-name for the empty set. 
\end{lemma}

\begin{proof}
	Let $\delta\colon \PP\to \lambda$ show that $\PP$ is explicitly $S$-closed. For $p\in \PP$ we let $\s(p) = 0^{\delta(p)+1}$, that is, a string of zeros of length $\delta(p)+1$. By well-ordering all (co-finality equivalence classes of) increasing sequences of conditions, we can make a choice of canonical exact upper bounds. 
\end{proof}



\subsection{Sparseness and iterations: the successor case} 


\begin{definition} \label{def:coherence}
Suppose that $\PP\compembeded \RR$, $S\subseteq \lambda$, that $\rho$ is an $S$-sparse $\PP$-name and that $\tau$ is an $S$-sparse $\RR$-name. We say that $\tau$ \emph{coheres} with~$\rho$ (and write $\rho\cohere \tau$) if:
\begin{orderedlist}
	\item If $p\in \dom \tau$ then $p\rest{\PP}\in \dom \rho$ and $\lth{\rho(p\rest{\PP})} = \lth{\tau(p)}$;\footnote{Note that we do not require that $\rho(p\rest{\PP}) = \tau(p)$.} 
	\item For an increasing sequence $\bar p$ from $\dom \tau$, if  $\cb_\tau(\bar p)$ is defined, then $\cb_\rho(\bar  p\rest{\PP})$ is defined and equals $\cb_\tau(\bar p)\rest{\PP}$. 
\end{orderedlist}
	\end{definition}
Note that if $\rho\cohere \tau$ and $\bar p$ is a sparse sequence for~$\tau$, then $\bar p\rest{\PP}$ is a sparse sequence for~$\rho$. 

\begin{definition} \label{def:one_step_sparse_name_iteration}
	Let~$S\subseteq \lambda$. Suppose that~$\PP$ is $\lllambda$-distributive; suppose that~$\rho$ is an $S$-sparse $\PP$-name; suppose that in~$V(\PP)$, $\QQ$ is a notion of forcing and~$\sigma$ is an $S\setminus W_{\rho}$-sparse $\QQ$-name. We define $\tau = \rho\vee \sigma$ as follows:
	\begin{itemize}
		\item $\dom \tau\subseteq \PP\starr\QQ$ is the collection of $(p,q)\in \PP\starr\QQ$ such that $p\in \dom{\rho}$, and, letting $\alpha = \lth{\rho(p)}$, there is some string $\pi\in 2^{<\lambda}$ (in~$V$) of length~$\alpha$ such that $p\force_{\PP} q\in \dom \sigma \andd \sigma(q) = \pi$. 
		\item For $(p,q)\in \dom \tau$ we define $\tau(p,q)$ to be the characteristic function of the union of $\rho(p)$ and $\sigma(q)$. That is, if $p\force \sigma(q)=\pi$ (where $\lth{\pi} = \alpha = \lth{\rho(p)}$) then we declare that $\lth{\tau(p,q)}= \alpha$, and for all $\beta< \alpha$, $\tau(p,q)(\beta)=1$ if and only if $\rho(p)(\beta)=1$ or $\pi(\beta)=1$. 
		\item If $(\bar p,\bar q)$ is an increasing sequence of conditions from $\dom \tau$, then we define $\cb_\tau(\bar p,\bar q) = (\cb_\rho(\bar p), \cb_\s(\bar q))$. That is, $\cb_\tau(\bar p,\bar q)$ is defined to be $(p^*,q^*)$ if $p^* = \cb_\rho(\bar p)$ is defined, and $p^*$ forces that $q^* = \cb_\s(\bar q)$.\footnote{The condition $\cb_\tau(\bar p,\bar q)$ is unique once we identify conditions $(p,q)$ and $(p,q')$ such that $p\force q = q'$.} 
	\end{itemize}
\end{definition}

\begin{lemma} \label{lem:iterating_one_step_with_a_sparse_forcing}
	Suppose that the hypotheses of \cref{def:one_step_sparse_name_iteration} hold: that $S\subseteq \lambda$, $\PP$ is $\lllambda$-distributive, $\rho$ is an $S$-sparse $\PP$-name, and in~$V(\PP)$, $\QQ$ is a notion of forcing and that~$\sigma$ is an $S\setminus W_{\rho}$-sparse $\QQ$-name. Suppose further that ~$S\cap \Cof(\aleph_0)$ is stationary. Then:
	\begin{itemize}
		\item $\rho\vee \sigma$ is an $S$-sparse $\PP\starr \QQ$-name;
		\item $\rho \cohere \rho\vee \sigma$; and
		\item in $V(\PP\starr\QQ)$, $W_{\rho\vee \sigma} = W_\rho\cup W_\sigma$. 
	\end{itemize}
\end{lemma}

\begin{proof}
	Let $\tau = \rho\vee \sigma$. 


	 
	Let us first show that $\dom \tau$ is dense in $\PP\starr\QQ$. Given $(p_{0},q_{0})\in \PP\starr \QQ$, since $S\cap \Cof(\aleph_0)$ is stationary, we find an increasing sequence of ordinals $\seq{\gamma_{n}}$ and a $\lambda$-filtration $\bar N$ such that $N_{\gamma_{n}}\cap \lambda = \gamma_{n}$ for all~$n$, $\PP,\QQ,\rho,\sigma,p_{0},q_{0}\in N_{{0}}$, and $\gamma_{\w} = \sup_{n} \gamma_{n}\in S$. We then define an increasing sequence of conditions $\seq{p_{n},q_{n}}\in \PP\starr \QQ$ \linebreak such that $(p_{n},q_{n})\in N_{\gamma_{n}}$, $p_{n}\in \dom \rho$, $\lth{\rho(p_{n})}\ge \gamma_{n-1}$, and $p_{n}$ forces that $q_n\in \dom \sigma$ and $\sigma(q_{n}) = \pi_{n}$ for some string $\pi_n\in V$ with $\gamma_{n-1}\le \lth{\pi_n}$. Here we use that~$\PP$ does not add sequences of ordinals of length $\lllambda$. Of course, since $(p_n,q_n)\in N_{\gamma_n}$, we have $\lth{\rho(p_n)},\lth{\pi_n}< \gamma_n$. 

	Let $\pi_\w = \bigcup_n \pi_n$. Since $\gamma_\w\in S$, $p_\w = \cb_\rho(\seq{p_{n}})$ is defined; and $p_{\w}$ forces that $\seq{q_{n}}$ is increasing in~$\QQ$ and that $\sigma(\bar q) = \pi_\w$ and so has length~$\gamma_\w$. Also, as $\rho(p_\w)(\gamma_\w)= 0$, $p_\w$ forces that $\gamma_\w\notin W_{\rho}$; so $p_\w$ forces that $q_\w=\cb_\s(\bar q)$ is defined. Note that $p_\w\force_\PP \sigma(q_\w) = \pi_\w\conc 0$.  Then $(p_{\w},q_{\w})\in \dom \tau$ and extends $(p_{0},q_{0})$. It is now not difficult to see that in $V(\PP\starr\QQ)$, $W_{\tau} = W_\rho\cup W_\sigma$. 

	\smallskip

	Next, we observe that if defined, $\cb_\tau(\bar p,\bar q)$ is in $\dom \tau$, and $\tau(\cb_\tau(\bar p,\bar q)) = \tau(\bar p,\bar q)\conc 0$; this is because $\rho(\cb_\rho(\bar p)) = \rho(\bar p)\conc 0$ and $\cb_\rho(\bar p)$ forces that $\s(\cb_\s(\bar q)) = \s(\bar q)\conc 0$. We also observe that $\cb_\tau$ is canonical. Finally, suppose that $(\bar p,\bar q)$ is sparse for~$\tau$, with $\alpha = \lth{\tau(\bar p,\bar q)}\in S$. Then $\bar p$ is sparse for~$\rho$, and $\lth{\rho(\bar p)}= \alpha$; let $p^* = \cb_\rho(\bar p)$. Then, $p^*$ forces that $\bar q$ is sparse for~$\s$, and that $ \lth{\s(\bar q)} = \alpha\in S\setminus W_\rho$; so $p^*$ forces that $\cb_\s(\bar q)$ is defined, whence $\cb_\tau(\bar p,\bar q)$ is defined. It follows that~$\tau$ is $S$-sparse, and that $\rho\cohere\tau$. 
\end{proof}





\Cref{lem:closed_implies_sparsely_closed} yields:

\begin{corollary} \label{lem:iterating_one_step_with_a_closed_forcing}
	Let~$S\subseteq \Cof_\lambda(\aleph_0)$ be stationary. Suppose that~$\rho$ is an $S$-sparse $\PP$-name; suppose that in~$V(\PP)$, $\QQ$ is a notion of forcing which is explicitly $S\setminus W_{\rho}$-closed. Then there is an $S$-sparse $\PP\starr\QQ$-name~$\tau$ for~$W_\rho$, which coheres with~$\rho$.
\end{corollary}


\subsection{Sparseness and iterations: the limit case} 

\begin{definition} \label{def:coherent_system_1}
	Let $\bar \PP = \seq{\PP_{i}}$ be an iteration. A \emph{coherent system of $S$-sparse names for~$\bar\PP$} is a sequence $\bar \tau = \seq{\tau_{i}}$ such that each~$\tau_{i}$ is an $S$-sparse $\PP_{i}$-name, and for $i<j < \lth{\bar \PP}$, $\tau_{i}\cohere \tau_{j}$. 
\end{definition}



\begin{definition} \label{def:theta_inverse_limit_of_coherent_system}
	Let $S\subseteq \lambda$. Suppose that $\theta<\lambda$ is regular, that $\bar \PP = \seq{\PP_{i}}_{i\le\theta}$ is an iteration with full support (a directed system with inverse limits), and that $\bar \tau = \seq{\tau_{i}}_{i<\theta}$ is a coherent system of $S$-sparse names for $\bar \PP\rest{\theta}$. We define $\tau_\theta = \bigvee \bar \tau$  as follows. First,
	\begin{itemize}
		\item $p\in \dom \tau_\theta$ if for all $i<\theta$, $p\rest{i}\in \dom \tau_{i}$.
	\end{itemize}
	Suppose that $p\in \dom \tau_\theta$. For $i<j<\theta$, because $\tau_i\cohere \tau_j$, we have $\lth{\tau_i(p\rest{i})} = \lth{\tau_j(p\rest{j})}$.
	\begin{itemize}
		\item For $p\in \dom \tau_\theta$, we let $\tau_\theta(p)$ be the characteristic function of the union of $\tau_i(p\rest{i})$. That is, if $\alpha = \lth{\tau_i(p\rest{i})}$ for all $i<\theta$, then we declare that $\lth{\tau_\theta(p)} = \alpha$, and for $\beta<\alpha$, $\tau_\theta(p)(\beta)=1$ if and only if there is some $i<\theta$ such that $\tau_i(p\rest{i})(\beta)=1$. 
		\item For an increasing sequence of conditions $\bar p$ from $\dom \tau_\theta$, we let $q = \cb_{\tau_\theta}(\bar p)$ if for all $i<\theta$, $q\rest{i} = \cb_{\tau_i}(\bar p\rest{i})$.
	\end{itemize}
\end{definition}

\begin{lemma} \label{lem:iterating_inverse_limit_case}
	Suppose that the hypotheses of \cref{def:theta_inverse_limit_of_coherent_system} hold: $S\subseteq \lambda$;  $\theta<\lambda$ is regular;  $\bar \PP = \seq{\PP_{i}}_{i\le\theta}$ is an iteration with full support, and $\bar \tau = \seq{\tau_{i}}_{i<\theta}$ is a coherent system of $S$-sparse names for $\bar \PP\rest{\theta}$. Suppose further that~$S$ is fat. 

	Then:
	\begin{itemize}
		\item  $\tau_\theta = \bigvee\bar \tau$ is an $S$-sparse $\PP_\theta$-name;
		\item for all $i<\theta$, $\tau_i \cohere \tau_\theta$; and
		\item  in $V(\PP_\theta)$, $W_{\tau_\theta} = \bigcup_{i<\theta} W_{\tau_i}$. 
	\end{itemize}
\end{lemma}

\begin{proof} 
	We show that~$\dom \tau_{\theta}$ is dense in~$\PP_{\theta}$. Let $r\in \PP_{\theta}$. Obtain a $\lambda$-filtration $\bar N$ and a set $D= \left\{ \gamma_{i} \,:\, i\le \theta  \right\}$ given by \cref{lem:approachable_ordinals}, with $S,r,\bar \PP, \bar \tau\in N_{0}$ (and $\theta\subset N_{0}$). We build a sequence $\bar p = \seq{p_{i}}_{i\le \theta}$ of conditions with the following properties:
	\begin{orderedlist}
		\item $p_{i}\in \PP_{i}$; 
		\item $p_{i}$ extends $r\rest i$; 
		\item If~$i$ is a successor then $p_i\in \dom \tau_i$ and $\lth{\tau_{i}(p_{i})}> \gamma_{i-1}$;
		\item \label{item:inverse_limit_construction:sparse_at_all_coordinates}
		for all $k<i$, $\seq{p_{j}\rest{k}}_{j\in (k,i)}$ is a sparse sequence for $\tau_{k}$.
	\end{orderedlist}
	Note that \ref{item:inverse_limit_construction:sparse_at_all_coordinates} implies that for all $j<i<\theta$, $p_i\rest{j}\in \dom \tau_j$. We do not however assume that for limit~$i$, $p_i\in \dom \tau_i$.

	We start with~$p_{0}$ being the empty condition. Given~$p_{i}$ we find an extension $p_{i+1}\in \dom \tau_{i+1}$ which also extends $r\rest{(i+1)}$. By extending, we can make $\lth{\tau_{i+1}(p_{i+1})}\ge \gamma_{i}$. 

	Now suppose that $i\le\theta$ is a limit ordinal. As we argued in the proof of \cref{prop:adding_shorter_sequences_with_fatness_and_explicit_closure}, for all $j<i$, $p_{j+1}\in N_{\gamma_j+1}$, and so $\lth{\tau_{j+1}(p_{j+1})}< \gamma_{j+1}$. Let $\alpha_i = \sup_{j<i}\gamma_j$; so $\alpha_i\in S$ (recall that the sequence $\seq{\gamma_i}$ need not be continuous). Now \ref{item:inverse_limit_construction:sparse_at_all_coordinates} above holds for~$i$ by induction; for all $k<i$, $|\tau_k(\seq{p_j\rest{k}}_{j\in (k,i)})| = \alpha_i$. It follows that for all~$k<i$, $q_k = q_k^i = \cb_{\tau_k}(\seq{p_j\rest{k}}_{j\in (k,i)})$ is defined. Further, if $k<k'<i$ then as $\cb_{\tau_k}$ is canonical, $q_k = \cb_{\tau_k}(\seq{p_j\rest{k}}_{j\in (k',i)})$. By coherence, $q_k = q_{k'}\rest{k}$. Since $\PP_{i}$ is the inverse limit of $\bar \PP\rest{i}$, we let $p_i\in \PP_i$ be the inverse limit of the sequence $\seq{q_k}_{k<i}$, that is, for all $k<i$, $q_k = p_i\rest{k}$. Note that \ref{item:inverse_limit_construction:sparse_at_all_coordinates} now holds for $i+1$. Also note that for all $k<i$, $q_k$ extends $p_k$ which extends $r\rest{k}$; it follows that $p_i$ extends $r\rest{i}$. At step~$\theta$ we get $p_\theta\in \dom \tau_\theta$ and extending~$r$.
	
	If $p,q\in \dom \tau_{\theta}$ and~$q$ extends~$p$, then for all $i<\theta$, $\tau_{i}(p\rest{i}) \preceq \tau_{i}(q\rest{i})$, so $\tau_{\theta}(p)\preceq \tau_{\theta}(q)$. In $V(\PP_{\theta})$, $W_{\tau_{\theta}} = \bigcup_{i} W_{\tau_{i}}$.

	\smallskip
	
	Let $\bar p$ be an increasing sequence from $\dom \tau_\theta$, and suppose that $p^* = \cb_{\tau_\theta}(\bar p)$ is defined. By definition, $p^*\in \dom \tau_\theta$. Let $\alpha = \lth{\tau_\theta(\bar p)}$. For all $i<\theta$, $\lth{\tau_i(\bar p\rest{i})} = \alpha$ and $\tau_i(p^*\rest{i})(\alpha)=0$, so by our definition, $\tau_\theta(p^*) = \tau_\theta(\bar p)\conc 0$. It is also easy to see that $\cb_\theta$ is canonical. 

	Suppose that $\bar p$ is a sparse sequence for~$\tau_{\theta}$, and that $\alpha = \lth{\tau_{\theta}(\bar p)}\in S$. For all $i<\theta$, $\bar p\rest{i}$ is sparse for~$\tau_i$, and so $q_i=\cb_{\tau_i}(\bar p\rest{i})$ is defined; and as above, for $i<i'<\theta$, $q_i = q_{i'}\rest{i}$. Then the inverse limit of $\seq{q_i}$ equals $\cb_{\tau_\theta}(\bar p)$.

	We conclude that $\tau_{\theta}$ is $S$-sparse and that $\tau_{i}\cohere \tau_{\theta}$ for all $i<\theta$. 
\end{proof}


For the next definition and lemma, note that if $\tau$ is an $S$-sparse $\PP$-name, and $A\subseteq \dom \tau$ is a dense final segment of $\dom \tau$, then $\tau\rest{A}$ is also an $S$-sparse $\PP$-name, and in $V(\PP)$, $W_\tau = W_{\tau\rest A}$.

\begin{definition} \label{def:theta_direct_limit_of_coherent_system}
	Let $S\subseteq \lambda$. Suppose that $\bar \PP = \seq{\PP_{i}}_{i\le\lambda}$ is an iteration with inverse limits below~$\lambda$ and a direct limit at~$\lambda$; suppose that $\bar \tau = \seq{\tau_{i}}_{i<\lambda}$ is a coherent system of $S$-sparse names for $\bar \PP\rest{\lambda}$. We define $\tau_\lambda = \bigvee \bar \tau$ as follows. First, 
	\begin{itemize}
		\item we let $\dom \tau_\lambda$ be the collection of $p\in \PP_\lambda$ such that for some limit $\alpha<\lambda$, $p\in \PP_\alpha$, and for all $\beta<\alpha$, $p\rest{\beta}\in \dom \tau_\beta$ and $\lth{\tau_\beta(p\rest{\beta})} = \alpha+1$.
	\end{itemize}
	For all $p\in \dom \tau_\lambda$ there is a unique~$\alpha$ witnessing this fact; $\alpha$ is determined by $\lth{\tau_\gamma(p\rest\gamma)}$ for any $\gamma<\alpha$.
	\begin{itemize}
		\item For $p\in \dom \tau_{\lambda}$, as witnessed by~$\alpha$, we declare that $\lth{\tau_{\lambda}(p)} = \alpha+1$ and that for $\beta \le \alpha$, $\tau_{\lambda}(p)(\beta)=1$ if and only if $\tau_{\gamma}(p\rest{\gamma})(\beta) = 1$ for some $\gamma<\beta$. 
	That is, $\tau_{\lambda}(p) = \DiagUnion_{\gamma<\alpha}\tau_{\gamma}(p\rest{\gamma})$. 
	\item Suppose that $\bar p$ is an increasing sequence from $\dom \tau_\lambda$. Let $\alpha = \lth{\tau_\lambda(\bar p)}$. We let $\cb_{\tau_\lambda}(\bar p) = p^*$ if $p^*\in \PP_\alpha$ and for all $\beta<\alpha$, $p^*\rest{\beta}$ is $\cb_{\tau_\beta}$ of a tail of $\bar p\rest{\beta}$. 
	\end{itemize}
\end{definition}

\begin{lemma} \label{lem:iterating_direct_limit_case}
	Suppose that the hypotheses of \cref{def:theta_direct_limit_of_coherent_system} hold: $S\subseteq \lambda$; $\bar \PP = \seq{\PP_{i}}_{i\le\lambda}$ is a $\lllambda$-support iteration; $\bar \tau$ is a coherent system of $S$-sparse names for $\bar \PP\rest{\lambda}$. Suppose further that $S\cap \Cof(\aleph_0)$ is stationary in~$\lambda$.

	Then: 
	\begin{itemize}
		\item $\tau_\lambda = \bigvee \bar \tau$ is an $S$-sparse $\PP_\lambda$-name;
		\item In $V(\PP_\lambda)$, $W_{\tau_\lambda} = \DiagUnion_{i<\lambda} W_{\tau_{i}}$; and
	\end{itemize}
	letting 
	\[
		C_{\alpha} = 
		\left\{ p\in \dom \tau_{\lambda} \,:\, \lth{\tau_{\lambda}(p)}>\alpha+1   \right\},
	\]
	\begin{itemize}
		\item for all $\alpha<\lambda$, $\tau_{\alpha} \cohere \tau_{\lambda}\rest{C_{\alpha}}$.
	\end{itemize}
\end{lemma}

Note that indeed each $C_\alpha$ is a dense final segment of $\dom \tau_\lambda$. 

\begin{proof}

	We show that $\dom \tau_{\lambda}$ is dense in~$\PP_{\lambda}$. Let $p_{0}\in \PP_{\lambda}$. Obtain an increasing sequence $\seq{\gamma_{n}}$ with $\gamma_{\w} = \sup_{n}\gamma_{n}\in S$, and models $N_{\gamma_{n}}$ with $N_{\gamma_{n}}\cap \lambda = \gamma_{n}$, such that all relevant information, including~$p_{0}$, is in $N_{\gamma_{0}}$. We then define an increasing sequence $\seq{p_{n}}$ with $p_{n}\in N_{\gamma_{n}}$ and for some $\alpha_{n}\in [\gamma_{n-1},\gamma_n)$,  $p_{n}\in \dom \tau_{\alpha_{n}}$,
	and $\lth{\tau_{\alpha_{n}}(p_{n})} \ge \alpha_{n}$. As usual, $\lth{\tau_{\alpha_n}(p_n)}< \gamma_n$. 	Now for all $\beta<\gamma_{\w}$, for all but finitely many~$n$ (say for all $n\ge n_{\beta}$), $p_{n}\rest{\beta}\in \dom \tau_{\beta}$; since $|\tau_\beta({\seq{p_n\rest{\beta}}_{n\ge n_\beta}})| = \gamma_\w\in S$, $q_\beta = \cb_{\tau_\beta}{\seq{p_n\rest{\beta}}_{n\ge n_\beta}}$ is defined. As above, this value does not change if we take a tail of the sequence, so for $\beta<\alpha<\gamma_\w$ we have $q_\beta = q_\alpha\rest{\beta}$; so the inverse limit of the sequence $\seq{q_\beta}$ is in $\dom \tau_\lambda$ (note that $\lth{\tau_\beta(q_\beta)} = \gamma_\w+1$).

	\smallskip
	
 Now for all $\alpha<\lambda$,   $p\rest{\alpha}\in \dom \tau_\alpha$ for all $p\in C_\alpha$. Further, by definition of $\tau_\lambda$, for $p\in C_\alpha$ we have $\lth{\tau_\lambda(p)} = \lth{\tau_\alpha(p\rest{\alpha})}$. 

	Suppose that $p,q\in \dom \tau_{\lambda}$ and that $q$ extends~$p$; say $\lth{\tau_{\lambda}(p)}= \alpha+1$ and $\lth{\tau_{\lambda}(q)}=\beta+1$. For any $\gamma< \alpha,\beta$, we have $p\rest\gamma \le q\rest{\gamma}$ and so $\tau_\gamma(p\rest{\gamma})\preceq \tau_\gamma(q\rest{\gamma})$; and $\alpha+1 = \lth{\tau_\gamma(p\rest{\gamma})}$, $\beta+1 = \lth{\tau_\gamma(q\rest{\gamma})}$. Hence $\alpha\le \beta$. 

	It is then not difficult to see that $\tau_{\lambda}(p)\preceq \tau_{\lambda}(q)$: to determine the value on $\gamma\le \alpha$ we note that for all $\delta<\gamma$, $\tau_{\delta}(p\rest{\delta})$ and $\tau_{\delta}(q\rest{\delta})$ agree on~$\gamma$. It follows that~$\tau_{\lambda}$ is an explicit~$\PP_\lambda$-name for $\DiagUnion_{i<\lambda} W_{\tau_{i}}$, and that each~$C_{\alpha}$ is a final segment of~$\dom \tau_{\lambda}$. 

	\smallskip
	
	Suppose that $\bar p$ is an increasing sequence from $\dom \tau_\lambda$, and let $\alpha = \lth{\tau_\lambda(\bar p)}$. Suppose that $p^* = \cb_{\tau_\lambda}(\bar p)$ is defined. Let $\beta<\alpha$. Then $\lth{\tau_\beta(p^*\rest{\beta})}= \alpha+1$ and $\tau_\beta(p^*\rest{\beta})(\alpha)=0$. By definition, $p^*\in \dom \tau_\lambda$, and $\tau_\lambda(p^*)(\alpha)=0$. Also, $\cb_{\tau_\lambda}$ is canonical.

	Suppose that $\bar p = \seq{p_i}_{i<i^*}$ is a sparse sequence for $\s_\lambda$. Let $\alpha = \lth{\tau_\lambda(\bar p)}$, and suppose that $\alpha\in S$.  For $\beta<\alpha$ let $i(\beta)$ be the least~$i$ such that $\lth{\tau_\lambda(p_i)} > \beta+1$. Then $\seq{p_i}_{i\in [i(\beta),i^*)}$ is a sequence from $C_\beta$, and $\seq{p_i\rest{\beta}}_{i\in [i(\beta),i^*)}$ is sparse for $\tau_\beta$. As usual, since $\PP_\alpha$ is an inverse limit we get $q\in \PP_\alpha$ such that for all $\beta<\alpha$, $q\rest{\beta} = \cb_{\tau_\beta}(\seq{p_i\rest{\beta}}_{i\in [i(\beta),i^*)})$; $q = \cb_{\tau_\lambda}(\bar p)$. 

	This argument shows that $\tau_{\lambda}$ is $S$-sparse and that $\tau_{\alpha}\cohere \tau_{\lambda}\rest{C_{\alpha}}$ for all~$\alpha$.
\end{proof}

Note that the sets $C_\alpha$ have a continuity property: suppose that $\bar p = \seq{p_i}$ is an increasing sequence (of limit length) from $\dom \tau_\lambda$ and $q = \cb_{\tau_\lambda}(\bar p)$ is defined. For $\alpha<\lambda$, if for all $i$, $p_i\notin C_\alpha$, then $q\notin C_\alpha$. Also note that for all $p\in \dom \tau_\lambda$, $p\in C_\alpha$ for fewer than~$\lambda$ many $\alpha$.


\subsection{Named $\lambda$-iterations} 
\label{sub:named_lambda_iterations}

\begin{definition} \label{def:named_iterations}
	We say that $\seq{\PP_{\zeta}}_{\zeta\le \xi}$, where $\xi\le \lambda^+$, is a \emph{named $\lambda$-iteration}, if there is a sequence $\seq{\QQ_\zeta,\s_\zeta}_{\zeta<\xi}$ such that:
	\begin{orderedlist}
		\item $\seq{\PP_\zeta,\QQ_\zeta}$ is a $\lllambda$-support iteration;
		\item For all $\zeta<\xi$, $\PP_\zeta$ is $\lllambda$-distributive;
		\item For all $\zeta<\xi$, $\sigma_\zeta\in V(\PP_\zeta)$ is an explicit $\QQ_\zeta$-name for a subset of~$\lambda$. 
	\end{orderedlist}
\end{definition}

Note that we are not requiring $\PP_\xi$ to be $\lllambda$-distributive; under further assumptions, this will follow, as we shall shortly see. We require that $\PP_\zeta$ be $\lllambda$-distributive for $\zeta<\lambda$ so that~$\lambda$ is regular in $V(\PP_\zeta)$, so that the notion of $\s_\zeta$ being an explicit $\QQ_\zeta$-name for a subset of~$\lambda$ (and later, a sparseness requirement) makes sense. When~$\xi$ is a limit ordinal, we also refer to the restriction $\PP\rest{\xi}$ as a named $\lambda$-iteration. 

\smallskip

If $\seq{\PP_\zeta}_{\zeta\le \xi}$ is a named $\lambda$-iteration, then we let, for all $\zeta\le\xi$, $\zeta<\lambda^+$, in $V(\PP_\zeta)$, $\uu_\zeta\in \PPP(\lambda)/\NS_\lambda$ be the least upper bound of 
\[
	\left\{ [W_{\sigma_\upsilon}]_{\NS_\lambda}  \,:\,  \upsilon<\zeta \right\}.
\]

\begin{definition} \label{def:sparse_iteration}
	Let $S\subseteq \lambda$. We say that a named $\lambda$-iteration $\seq{\PP_\zeta}_{\zeta\le \xi}$ is \emph{$S$-sparse} if for every $\zeta<\xi$, in $V(\PP_\zeta)$, $\sigma_\zeta$ is $S\setminus \uu_\zeta$-sparse.\footnote{We mean that it is $S\setminus U$-sparse for some~$U$ such that $[U]_{\NS_\lambda} = \uu_\zeta$. By \cref{lem:sparseness_invariant_under_clubs}, up to a slight modification of $\sigma_\zeta$, the choice of~$U$ does not matter.} 
\end{definition}

Our aim is to show that when~$S$ is fat, if $\seq{\PP_\zeta}$ is an $S$-sparse iteration, then for each $\zeta\le \xi$, $\zeta<\lambda^+$, there is an $S$-sparse $\PP_\zeta$-name $\tau_\zeta$ such that (in $V(\PP_\zeta)$) $[W_{\tau_\zeta}]_{\NS_\lambda} = \uu_\zeta$. However, because we may have $\xi>\lambda$, we will not be able to get a coherent sequence $\bar \tau$ of names. We need a weak notion of coherence:

\begin{definition} \label{def:weak_coherent_sequence}
	 Suppose that $\bar \PP = \seq{\PP_\zeta}_{\zeta<\xi}$ is an iteration, and that $\bar \tau = \seq{\tau_\zeta}_{\zeta<\xi}$ is a sequence such that for all $\zeta<\xi$, $\tau_\zeta$ is an explicit $\PP_\zeta$-name for a subset of~$\lambda$. We say that $\bar \tau$ is \emph{weakly coherent} if for $\upsilon\le \zeta<\xi$ there are sets $A^\zeta_\upsilon$ satisfying:
	\begin{sublemma}
			\item Each $A^\zeta_\upsilon$ is a dense final segment of $\dom \tau_\zeta$; 
			\item For all $\upsilon\le \zeta<\xi$, $\tau_\upsilon \cohere \tau_\zeta\rest{A^\zeta_\upsilon}$; 
			\item $A^\zeta_\zeta = \dom \tau_\zeta$;
			\item \label{item:weak_coherence:first_restriction}
			For $\alpha\le \beta\le \gamma< \xi$,  $p\in A^\gamma_\alpha\cap A^\gamma_\beta$ implies $p\rest{\beta}\in A^\beta_\alpha$;
			\item \label{item:weak_coherence:second_restriction}
			For $\alpha\le \beta\le \gamma< \xi$,  $p\in A^\gamma_\beta$ and $p\rest{\beta}\in A^\beta_\alpha$ implies $p\in A^\gamma_\alpha$;
			\item \label{item:weak_coherence:number_bound}
			For all $\zeta<\xi$ and every $p\in \dom \tau_\zeta$, there are $\lllambda$ many $\upsilon<\zeta$ such that $p\in A^\zeta_\upsilon$;
			\item \label{item:weak_coherence:continuity}
			For $\upsilon<\zeta$, if $\bar p = \seq{p_i}$ is an increasing sequence from $\dom \tau_\zeta$,  $p = \cb_{\tau_\zeta}(\bar p)$ is defined, and for all~$i$, $p_i\notin A^\zeta_\upsilon$, then $p\notin A^\zeta_\upsilon$. 
	\end{sublemma}
\end{definition}

We will show that we can construct a weakly coherent sequence of names as required. We will use more properties of the sequence, which we incorporate into the following definition. 

\begin{definition} \label{def:associated_tau_sequence_regularity_properties:up_to_lambda_plus}
	Let $\seq{\PP_\zeta}_{\zeta<\xi}$ be an $S$-sparse iteration. An \emph{associated sequence} is a sequence $\bar \tau = \seq{\tau_\zeta}_{\zeta\in [1,\xi)}$ such that:
	\begin{equivalent}
		\item Each $\tau_\zeta$ is an $S$-sparse $\PP_\zeta$-name; 
		\item \label{item:lambda_plus_iteration:U_s_correct}
		In $V(\PP_\zeta)$, $[W_{\tau_\zeta}]_{\NS_\lambda} = \uu_\zeta$; 
		\item  \label{item:lambda_plus_iteration:weak_coherence}
		The sequence $\bar \tau$ is weakly coherent;
		\item \label{item:lambda_plus_iteration:uniform_length_of_conditions}
		For $p\in \dom \tau_\zeta$, for all $\upsilon< \zeta$, $p\in A^\zeta_{\upsilon+1}$ if and only if $\upsilon\in \supp(p)$, in which case there is a string $\pi_\upsilon\in 2^{<\lambda}$ of length $\lth{\tau_\zeta(p)}$ such that $p\rest{\upsilon}\force_{\PP_\upsilon} p(\upsilon)\in \dom \s_\zeta \andd \sigma_\upsilon(p(\upsilon)) = \pi_\upsilon$;
	\end{equivalent}
	We write $\s_\upsilon(p(\upsilon))$ for the string $\pi_\upsilon$;
	\begin{equivalent}[resume]
		\item \label{item:lambda_plus_iteration:subset_relations_from_a_point}
		If $p,q\in \dom \tau_\zeta$, $q$ extends~$p$, and $\upsilon\in \supp(p)$, then $\s_\upsilon(q(\upsilon))\setminus \lth{\tau_\zeta(p)}\subseteq \tau_\zeta(q)$;\footnote{Again, recall that we identify sets and characteristic functions; so this means: for all $\alpha$ with $|\tau_\zeta(p)|\le \alpha < |\tau_\zeta(q)|$, if $\sigma_\upsilon(q(\upsilon))(\alpha)=1$ then $\tau_\zeta(q)(\alpha)=1$.}

		\item \label{item:lambda_plus_iteration:zero_in_sigmas_implies_in_tau}
		If $p\in \dom \tau_{\zeta}$,  $\alpha<\lth{\tau_\zeta(p)}$, and for all $\upsilon\in \supp(p)$, $\s_\upsilon(p(\upsilon))(\alpha)=0$, then $\tau_\zeta(p)(\alpha)=0$;

		\item \label{item:lambda_plus_iteration:continuity_of_support_in_canonical_bounds}
		If $\bar p = \seq{p_i}$ is an increasing sequence from $\dom \tau_\zeta$ and $p = \cb_{\tau_\zeta}(\bar p)$ is defined, then $\supp(p) = \bigcup_i \supp(p_i)$;

		\item \label{item:lambda_plus_iteration:taking_canonical_bounds_at_every_location}
		If $\bar p = \seq{p_i}$ is an increasing sequence from $\dom \tau_\zeta$ and $p = \cb_{\tau_\zeta}(\bar p)$ is defined, then for all $\upsilon\in \supp(p)$, $p\rest{\upsilon}$ forces that $p(\upsilon)$ is $\cb_{\s_\upsilon}$ of a tail of $\seq{p_i(\upsilon)}$. 
	\end{equivalent}
\end{definition}

\begin{notation}
	If $\bar \PP$ is an $S$-sparse iteration, and $\bar \tau$ is an associated sequence, then we write $W_\zeta$ for $W_{\s_\zeta}$ and $U_\zeta$ for $W_{\tau_\zeta}$. 
\end{notation}

\begin{proposition} 
\label{prop:iteration_up_to_lambda_plus:existence_of_associated_sequence:inductive_step}
Let $S\subseteq \lambda$ be fat. Suppose that $\seq{\PP_\zeta}_{\zeta\le \xi}$ is an $S$-sparse named $\lambda$-iteration, with $\xi<\lambda^+$, and that $\bar \tau = \seq{\tau_\zeta}_{\zeta\in [1,\xi)}$ is an associated sequence for $\bar \PP\rest{\xi}$. 

Then there is some $\tau_\xi$ such that $\bar \tau\conc \tau_\xi$ is an associated sequence for $\bar \PP$. 
\end{proposition}

In particular, $\PP_\xi$ is $\lllambda$-distributive.

\begin{corollary} \label{cor:iteration_up_to_lambda_plus:existence_of_associated_sequence}
	If~$S$ is fat, $\xi\le \lambda^+$, and $\bar \PP = \seq{\PP_\zeta}_{\zeta<\xi}$ is an $S$-sparse iteration, then $\bar \PP$ has an associated sequence~$\bar \tau$. 
\end{corollary}

\begin{proof}[Proof of \cref{prop:iteration_up_to_lambda_plus:existence_of_associated_sequence:inductive_step}]
	The definition of $\tau_\xi$ is of course by cases. 

	\medskip
	
	\noindent\underline{Case I: $\xi = 1$}. We let $\tau_1 = \s_0$ (recall that $\PP_1 = \QQ_0$). 

	\medskip
	
	\noindent\underline{Case II: $\xi = \vartheta+1$, $\vartheta>0$}. 
	We apply \cref{lem:iterating_one_step_with_a_sparse_forcing} to $\PP = \PP_\vartheta$, $\rho = \tau_\vartheta$, $\QQ = \QQ_\vartheta$, and~$\sigma$ being some mild variation of $\sigma_\vartheta$ which is sparsely $S\setminus U_\vartheta$-closed; let $\tau_{\xi}$ be the~$\tau$ obtained. For $\zeta\le \vartheta$, we let $A^\xi_\zeta =  \big\{ p\in \dom \tau_\xi \,:\,  p\rest{\vartheta} \in A^\vartheta_\zeta \big\}$; note that $A^\xi_\vartheta  = \dom \tau_\xi$. 

	For \ref{item:lambda_plus_iteration:U_s_correct} of \cref{def:associated_tau_sequence_regularity_properties:up_to_lambda_plus}, note that $\uu_\xi = \uu_\vartheta\vee [W_{\vartheta}]_{\NS_\lambda}$, and $U_\xi = U_\vartheta\cup W_\vartheta$. \ref{item:lambda_plus_iteration:weak_coherence} is not difficult. For \ref{item:lambda_plus_iteration:uniform_length_of_conditions}, note that for all $p\in \dom \tau_\xi$, $\vartheta\in \supp(p)$ and $p\rest{\vartheta}\force p(\vartheta) = \pi$ for some~$\pi$ of length $\lth{\tau_\xi(p)} = \lth{\tau_\vartheta(p\rest{\vartheta})}$. \ref{item:lambda_plus_iteration:subset_relations_from_a_point} and \ref{item:lambda_plus_iteration:zero_in_sigmas_implies_in_tau} follow from $\tau_\xi(q) = \tau_\vartheta(p\rest{\vartheta})\cup \s_\vartheta(p(\vartheta))$. \ref{item:lambda_plus_iteration:continuity_of_support_in_canonical_bounds} follows from the fact that $\cb_{\tau_\vartheta}(\bar p\rest{\vartheta}) = p\rest{\vartheta}$, and that $\vartheta\in \supp(q)$ for all $q\in \dom \tau_\xi$. \ref{item:lambda_plus_iteration:taking_canonical_bounds_at_every_location} follows similarly, noting that \cref{def:one_step_sparse_name_iteration} implies that it holds for $\upsilon = \vartheta$. 

	\medskip
	
	\noindent\underline{Case III: $\xi$ is a limit ordinal and $\cf(\xi)<\lambda$}. 
	For all $\zeta<\xi$, $\PP_\zeta$ is $\lllambda$-distributive. Then for all $J\subseteq \zeta$, if $|J|<\lambda$ then 
	\[
		A^\zeta_J = \bigcap_{\upsilon\in J} A^\zeta_\upsilon
	\]
	is a dense final segment of $\dom \tau_\zeta$; for all $\upsilon\in J$, $\tau_\upsilon \cohere \tau_\zeta\rest{A^\zeta_J}$. Further, for $J\subseteq \xi$ with $|J|<\lambda$, the sequence 
	\[
		\seq{\tau_\zeta \rest{A^\zeta_{J\cap \zeta}}}_{\zeta\in J}
	\]
	is coherent; this follows from condition~\ref{item:weak_coherence:first_restriction} of \cref{def:weak_coherent_sequence}. 

	We fix a closed, unbounded $J\subseteq \zeta$ of order-type $\theta = \cf(\zeta)$, and apply \cref{lem:iterating_inverse_limit_case} to the sequence $\smallseq{\tau_\zeta \rest{A^\zeta_{J\cap \zeta}}}_{\zeta\in J}$; notice that for limit points $\delta$ of~$J$, as $\cf(\delta)<\lambda$, $\PP_\delta$ is the inverse limit of $\seq{\PP_{\upsilon}}_{\upsilon\in J\cap \delta}$. We let $\tau_\xi$ be the~$\tau$ obtained. So for all $\delta\in J$, $\tau_\delta\rest{A^\delta_{J\cap \delta}}\cohere \tau_\xi$. 

	For \ref{item:lambda_plus_iteration:U_s_correct}, we use $U_{\xi} = \bigcup_{\delta\in J} U_\delta$, and the fact that $J$ is cofinal in~$\xi$. For \ref{item:lambda_plus_iteration:weak_coherence}, we notice that for all $\gamma<\xi$, for $p\in \dom \tau_\xi$, there is some $\delta\in J$, $\delta\ge \gamma$ such that $p\rest{\delta}\in A^\delta_\gamma$ if and only if for all $\delta\in J$ such that $\delta\ge \gamma$, $p\rest{\delta}\in A^\delta_\gamma$; we use either~\ref{item:weak_coherence:first_restriction} or~\ref{item:weak_coherence:second_restriction} of \cref{def:weak_coherent_sequence}; in either case we note that for all $\delta<\epsilon$ from~$J$, for all $p\in \dom \tau_\xi$, $p\rest{\epsilon}\in A^\epsilon_\delta$. We let $p\in A^\xi_\gamma$ if these equivalent conditions hold. Note that $A^\xi_\delta = \dom \tau_\xi$ for all $\delta\in J$. For \ref{item:weak_coherence:number_bound}, we use that $\bar \tau$ is an associated sequence, and $|J|<\lambda$. \ref{item:lambda_plus_iteration:subset_relations_from_a_point} and \ref{item:lambda_plus_iteration:zero_in_sigmas_implies_in_tau} of \cref{def:associated_tau_sequence_regularity_properties:up_to_lambda_plus} follow from  $\tau_\xi(q) = \bigcup_{\delta\in J} \tau_\delta(p\rest{\delta})$. \ref{item:lambda_plus_iteration:continuity_of_support_in_canonical_bounds} follows from the fact that if $p = \cb_{\tau_\xi}(\bar p)$ then for all $\delta\in J$, $p\rest{\delta} = \cb_{\tau_\delta}(\bar p\rest{\delta})$, and of course $\supp(q) = \bigcup_{\delta\in J}\supp(q\rest{\delta})$ for all $q\in \PP_\xi$. \ref{item:lambda_plus_iteration:taking_canonical_bounds_at_every_location} follows from \ref{item:lambda_plus_iteration:continuity_of_support_in_canonical_bounds} and our assumptions on~$\bar \tau$.

	\medskip
	
	\noindent\underline{Case IV: $\cf(\xi) = \lambda$}. 
	The construction is similar to case~(III). We fix $J \subseteq \xi$ closed and unbounded of order-type~$\lambda$; let $\seq{\delta_i}_{i<\lambda}$ be the increasing enumeration of~$J$. We again have that $\smallseq{\tau_{\delta}\rest{A^\delta_{J\cap \delta}}}_{\delta\in J}$ is a coherent system, with $\PP_\delta$ being an inverse limit of $\seq{\PP_\upsilon}_{\upsilon\in J\cap \delta}$ for all limit points~$\delta$ of~$J$, whereas this time, $\PP_\xi$ is the direct limit of $\seq{\PP_\delta}_{\delta\in J}$. So this time we apply \cref{lem:iterating_direct_limit_case} to get~$\tau_\xi$. In this case, for $\delta = \delta_i\in J$ we obtain the sets 
	\[
		C_\delta = \left\{ p\in \dom \tau_\xi \,:\, \lth{\tau_\xi(p)}> i+1  \right\}
	\]
	So if $\lth{\tau_\xi(p)}= j+1$ then $p\in C_{\delta_i}\Iff i<j$, and $\tau_{\delta}\rest{A^\delta_{J\cap \delta}}\cohere \tau_\xi\rest C_\delta$ for all $\delta\in J$. 

	\smallskip
	
	For \ref{item:lambda_plus_iteration:U_s_correct}, we use the fact that $U_\xi = \DiagUnion U_{\delta_i}$ so $\uu_\xi = \sup_{\delta\in J} \uu_\delta$. Toward \ref{item:lambda_plus_iteration:weak_coherence}, let $\gamma<\xi$, let $p\in \dom \tau_\xi$, and let $\lth{\tau_\xi(p)}= j+1$. If $\gamma<\delta_j$, then the following are equivalent: $p\rest{\delta}\in A^\delta_\gamma$ for some $\delta\in J\cap[\gamma,\delta_j)$; and  $p\rest{\delta}\in A^\delta_\gamma$ for all $\delta\in J\cap[\gamma,\delta_j)$. We let $p\in A^\xi_\gamma$ if these equivalent conditions hold. If $\gamma\ge \delta_j$ then $p\notin A^\xi_\gamma$. The argument is then the same as in case~(III), restricting to $\delta< \delta_j$; \ref{item:weak_coherence:number_bound} follows from $j<\lambda$. For \ref{item:weak_coherence:continuity}, letting $j+1 = \lth{\tau_\xi(p)}$, assuming that $\upsilon<\delta_j$, we fix some $\delta\in J\cap [\gamma,\delta_j)$; for a tail of~$\bar p$ we have $p_i\in C_\delta$; we use the fact that $\cb_{\tau_\xi}(\bar p)$ does not change when restricting to this tail. 

	For \ref{item:lambda_plus_iteration:uniform_length_of_conditions}, if $\lth{\tau_\xi(p)}= j+1$ then by the construction in \cref{lem:iterating_direct_limit_case}, $p\in \PP_{\delta_j}$, in other words $\supp(p)\subseteq \delta_j$. Also $p\in A^\xi_\upsilon$ implies $\upsilon < \delta_j$. For $\upsilon<\delta_j$, we find $\delta\in J\cap (\upsilon,\delta_j)$, and use \ref{item:lambda_plus_iteration:uniform_length_of_conditions} for $\tau_\delta$. 

	For \ref{item:lambda_plus_iteration:subset_relations_from_a_point} and \ref{item:lambda_plus_iteration:zero_in_sigmas_implies_in_tau} we use the definition of $\tau_\xi(p)$: if $\lth{\tau_\xi(p)}= j+1$ then $\tau_\xi(p)(\beta)=1$ if and only if $\tau_{\delta_i}(p\rest{\delta_i})(\beta)=1$ for some $i<\beta$. If $i<\beta$ and for all $\upsilon\in \supp(p\rest{\delta_i})$, $\s_\upsilon(p(\upsilon))(\beta)=0$, then $\tau_{\delta_i}(p\rest{\delta_i})(\beta)=0$. On the other hand, suppose that $q$ extends~$p$, $\nu\in \supp(p)$, that $\beta\in [\lth{\tau_\xi(p)}, \lth{\tau_\xi(q)})$, and that $\s_\upsilon(q)(\beta)=1$. Let $j+1 = \lth{\tau_{\xi}(p)}$. Since $\upsilon\in \supp(p)$, $\upsilon < \delta_j$. Choose $\delta\in J\cap [\upsilon,\delta_j)$. Then $p\rest{\delta}, q\rest{\delta}\in \dom \tau_{\delta}$, so by assumption on~$\bar \tau$, $\tau_\delta(q\rest{\delta})(\beta)=1$. Since $\beta\ge \delta_j > \delta$, by definition, $\tau_\xi(q)(\beta)=1$. 

	Finally, for \ref{item:lambda_plus_iteration:continuity_of_support_in_canonical_bounds}, let $p = \cb_{\tau_\xi}(\bar p)$; let $j+1 = \lth{\tau_\xi(p)}$. Then by construction, $\supp(p)= \bigcup_{\delta\in J\cap \delta_j} \supp (\cb_{\tau_\delta}(\bar p\rest{\delta}))$ (where we actually take a tail of $\bar p\rest{\delta}$). Again \ref{item:lambda_plus_iteration:taking_canonical_bounds_at_every_location} follows. 
\end{proof}

\subsection{Some more on named iterations} 
\label{sub:some_more_on_named_iterations}

\begin{lemma} \label{lem:preservation_of_strategic_S_closure}
	Let $\xi$ be a limit ordinal with $\cf(\xi)\ge \lambda$, suppose that $\seq{\PP_\zeta}_{\zeta<\xi}$ is a $\lllambda$-support iteration, that each $\PP_\zeta$ is $\lllambda$-distributive, and that $\bar \tau$ is a weakly coherent, $S$-sparse sequence for~$\bar \PP$. Then~$\PP_\xi$ is strategically $S$-closed. 
\end{lemma}

As discussed above, by the generalisation of \cref{prop:adding_shorter_sequences_with_fatness_and_explicit_closure} to strategic closure, this implies that if~$S$ is fat then $\PP_\xi$ is $\lllambda$-distributive. 

\begin{proof}
	The strategy for the completeness player is to play sequences $\bar p = \seq{p_i}$, where $p_i\in \PP_{\zeta_i}$, with $\seq{\zeta_i}$ increasing and continuous, such that:
	\begin{itemize}
		\item If the completeness player plays $p_{i+1}$, then $p_{i+1}\in \bigcap_{j\le i} A^{\zeta_{i+1}}_{\zeta_j}$; the ordinal played is $\lth{\tau_{\zeta_{i+1}}(p_{i+1})}$. 
		\item For limit~$i$, for all $k<i$, $p_i\rest{\zeta_k} = \cb_{\tau_{\zeta_k}} (\smallseq{p_j\rest{\zeta_k}}_{j\in (k,i)})$. \qedhere
	\end{itemize}
\end{proof}


\begin{lemma} \label{lem:the_iteration_lemma}
	Suppose that $\seq{\PP_\zeta}_{\zeta<\xi}$ is an $S$-sparse iteration, and that~$S$ is fat. For all $\upsilon  <\xi$, in $V(\PP_\upsilon)$, the iteration $\seq{\PP_\zeta/\PP_\upsilon}_{\zeta\in [\upsilon,\xi)}$ is $S\setminus U_\upsilon$-sparse (witnessed by $\seq{\QQ_\zeta,\s_\zeta}_{\zeta\in [\upsilon,\xi)}$). 
\end{lemma}

\begin{proof}
	All we need to observe is that as $\PP_\upsilon$ is $\lllambda$-distributive, the quotient iteration is $\lllambda$-support. 
\end{proof}

\begin{lemma} \label{lem:extending_compatibles_and_incompatibles}
	Suppose that $\seq{\PP_\zeta}_{\zeta\le \xi}$ is an $S$-sparse iteration, that $\xi<\lambda^+$, that~$S$ is fat, and let $\bar \tau$ be an associated sequence. 

	Let $\upsilon <\xi$. Suppose that $\kappa<\lambda$ and that $\left\{ p_i \,:\,  i<\kappa \right\}\subseteq \PP_\xi$ is a family of conditions such that $p_i\rest{\upsilon} = p_j\rest{\upsilon}$ for all $i,j<\kappa$. Then there is a collection of conditions $\left\{ q_i \,:\,  i<\kappa \right\}$ such that:
	\begin{orderedlist}
		\item Each~$q_i$ extends~$p_i$;
		\item For all~$i$, $q_i\in \dom \tau_\xi$, indeed $q_i\in A^\xi_{\upsilon}$
		\item For $i<j<\kappa$, $q_i\rest{\upsilon}= q_j\rest{\upsilon}$ 
	\end{orderedlist}
\end{lemma}

Note that it follows that for all $i<j<\kappa$, $\lth{\tau_\xi(q_i)} = \lth{\tau_\xi(q_j)}$.

\begin{proof}
	By induction on~$\kappa$. First we consider $\kappa=k$ finite. We define a sequence $\seq{r_n}$ of conditions such that $\seq{r_n\rest{\upsilon}}$ is increasing, and as usual $\sup_n \lth{\tau_\upsilon(r_n\rest{\upsilon})}\in S$. We start with $r_i = p_i$ for $i<k$; then at step $n = i\mod k$, we find an extension of $r_{n-1}\rest{\upsilon}\vee r_{n-k}$ in $A^\xi_\upsilon$. At the limit, for $i<k$ we let $q_i = \cb_{\tau_\xi}(\seq{r_{kn+i}}_{n<\w})$. The sequences $\seq{r_{kn+i}\rest{\upsilon}}$ are co-final and so have the same canonical upper bound. 

	\smallskip
	
	For infinite~$\kappa<\lambda$, let $\theta = \cf(\kappa)$, and let $\seq{\alpha_i}_{i\le \theta}$ be a continuous increasing sequence converging to~$\kappa$. We build an array of conditions $\smallseq{r^j_i}_{i\le \theta, j<\alpha_i}$ such that:
	\begin{itemize}
		\item For each $j<\kappa$, say $j\in [\alpha_{k},\alpha_{k+1})$, the sequence $\smallseq{r^j_i}_{i
		>k}$ is a sparse sequence for $\tau_\xi$;
		\item Each $r^j_i$ is in $A^\xi_\upsilon$; 
		\item For $j<j'<\alpha_i$, $r^j_i\rest\upsilon = r^{j'}_i\rest{\upsilon}$.
		\item For each $i<\theta$ and $j<\alpha_i$, $r^j_i$ extends $p_j
		$. 
	\end{itemize}
	Let~$s_i$ be the common value of $r^j_i\rest{\upsilon}$; it follows that $\seq{s_i}$ is sparse for $\tau_\upsilon$. At step $k+1<\theta$ we apply the inductive hypothesis to the collection of conditions
	\[
		\big\{ r^j_k \,:\,  j<\alpha_k \big\}  \cup \big\{  s_k\vee p_j \,:\, j\in [\alpha_k,\alpha_{k+1})    \big\}.
	\]
	As usual everything is happening within a filtration $\bar N$ of length $\theta+1$ with limit points in~$S$, so at limit steps we can take canonical upper bounds. We let $q_j = r^j_\theta$. 
\end{proof}

\section{Tasks and the task iteration} 
\label{sec:tasks_and_the_task_iteration}


\subsection{Tasks} 
\label{sub:tasks}

\begin{notation}
	Let $\FF_\lambda$ be the collection of $\lllambda$-distributive notions of forcing. Each $\PP\in \FF_\lambda$ preserves $\HH_\lambda$.
\end{notation}

	Note that since we assume that $\GCH$ holds below~$\lambda$, we may identify $\PPP(\lambda)$ with $\PPP(\HH_{\lambda})$.

\subsubsection{Generically uniform definitions} 

Let $O\subseteq \lambda$ and let $\vphi$ be a first-order formula. For every $A\in \PPP(\lambda)$ we define 
\[
	\BB(A) = \BB_{O,\vphi}(A) = \left\{ x\in \HH_\lambda 	\,:\, (\HH_\lambda;O,A)\models \vphi(x)   \right\}.
\]
We call this a uniform definition of $\BB(A)$ from~$A$; $O$ can be considered as an ``oracle'', but can also incorporate any parameters from $\HH_\lambda$ if we so choose. 

Further, for every $\PP\in \FF_\lambda$ and every $A\in \PPP(\lambda)(\PP)$, that is, any $A\in V(\PP)$ such that $\force_{\PP} A\subseteq \lambda$, we let $\BB(A)$ be a $\PP$-name denoting the result of this definition in $V(\PP)$. Thus, $O$ and~$\vphi$ give us what we call a \emph{generically uniform definition} of subsets of $\HH_\lambda$. 

We remark that if $\PP\compembeded \QQ$ and $\QQ\in \FF_\lambda$ then as $\HH_\lambda$ is the same in $V(\PP)$ and in $V(\QQ)$, we can naturally consider $\PPP(\lambda)(\PP)$ as a subset of $\PPP(\lambda)(\QQ)$. Since~$\vphi$ is first-order, the interpretation of $\BB(A)$ in $V(\PP)$ and $V(\QQ)$ is the same. 



\subsubsection{Limited genericity} 

\begin{definition}
	Suppose that $\PP\subseteq \HH_\lambda$ is a notion of forcing and that $O\in \PPP(\lambda)$. A filter $G\subseteq \PP$ is \emph{$O$-generic} if it meets every dense subset of~$\PP$ which is first-order definable in the structure $(\HH_\lambda;\PP,O)$ (we allow parameters from~$\HH_\lambda$).
\end{definition}

Note that if $\PP$ is $\lambda$-strategically closed then for any~$O$ there is an $O$-generic $G\subseteq \PP$ in~$V$.\footnote{Strategic $S$-closure for a fat~$S$ does not seem to suffice.}

If $\s$ is an explicit $\PP$-name for a subset of~$\lambda$, then for any $\s$-generic filter $G\subseteq \PP$, $W_\s[G]$ is well-defined and is an element of $2^\lambda$.

\subsubsection{Tasks} 
\label{ssub:tasks}


\begin{definition} \label{def:task}
	A \emph{$\lambda$-task} $\task$ consists of three generically uniform definitions $\QQ^\task,\s^\task$ and $\SS^\task$ (say with oracle~$O = O^\task$) such that for all $\PP\in \FF_\lambda$ and all $C\in \PPP(\lambda)(\PP)$, in $V(\PP)$, 
	\begin{enumerate}
		\item[(A)] $\QQ^\task(C)$ is a notion of forcing, and~$\s^\task(C)$ is a $\lambda$-sparse $\QQ^{\task}(C) $-name for a subset of~$C$.\footnote{This means that $\force_{\QQ^\task(C)} W_{\s^\task(C)}\subseteq C$.} We let $W^\task(C) = W_{\s^\task(C)}$. 
		\item[(B)] For any $G\subseteq \QQ^\task(C)$ which is $(O,C)$-generic, for any $\PP'\in \FF_\lambda$ with $\PP\compembeded \PP'$, for any $A\in \PPP(\lambda)(\PP')$, in $V(\PP')$,  
		 $\SS^\task(C,G,A)$ is a notion of forcing which is explicitly closed outside $W^\task(C)[G]$.\footnote{Recall that this means that it is explicitly $\lambda\setminus W^\task(C)[G]$-closed.}
	\end{enumerate}
\end{definition}

Note that explicit closure or sparse closure of a notion of forcing does not depend on the universe we work in, as it only involves sequences of length $<\!\lambda$.

\subsubsection{Example: Exact diamonds} 
\label{ssub:example_exact_diamonds}

\begin{definition}
	Let $\mu<\lambda$ be a cardinal (possibly finite, but $\ge 2$), and let $S\subseteq\lambda$ be stationary. A $\mu$-sequence on~$S$ is a sequence $\bar F = \seq{F_\alpha}_{\alpha\in S}$ such that for each $\alpha\in S$, $F_\alpha\subseteq \PPP(\alpha)$ and $|F_\alpha|= \mu$. 

A $\mu$-sequence $\bar F$ \emph{guesses} a set $X\subseteq \lambda$ if for stationarily many $\alpha\in S$, $X\cap \alpha\in F_\alpha$. 

We say that $\bar F$ is a \emph{$\mu$-diamond sequence} if it guesses every $X\subseteq \lambda$. 

We say that it is a \emph{$\mu$-almost diamond sequence} if it guesses all subsets of~$\lambda$ except possibly for a collection of at most~$\lambda$ many subsets~$X$. 

We say that $\bar F$ is a \emph{$\mu$-exact diamond sequence} if it is a $\mu$-diamond sequence, and 
whenever we choose $x_\alpha\in F_\alpha$ for all $\alpha\in S$, the sequence $\seq{F_\alpha\setminus \{x_\alpha\}}_{\alpha\in S}$ is not a $\mu$-almost diamond sequence (let alone a $\mu$-diamond sequence).
\end{definition}

We define the \emph{$\mu$-exact diamond task} $\task = \task_{\Diamond}(\mu)$. The forcing $\QQ^\task(C)$ adds a $\mu$-diamond sequence on a subset of~$C$; the subsequent forcings~$\SS^\task$ ensure that it is exact.

\begin{orderedlist}
	\item For $C\subseteq \lambda$, $p\in \QQ^\task(C)$ if $p = (\s(p),\bar F(p))$ where:
	\begin{itemize}
		\item $\s(p)\in 2^{<\lambda}$ and $\s(p)\subseteq C$;\footnote{As usual, identifying the characteristic function with its underlying set.} 
		\item $\bar F(p) = \seq{F_\alpha(p)}_{\alpha\in \s(p)}$ with $F_\alpha(p)\subseteq \PPP(\alpha)$ and $|F_\alpha(p)|= \mu$; 
		\item extension is by extending the sequences in both coordinates. 
	\end{itemize}
	$\s^\task(C) = \s$; i.e., the function sending $p$ to $\s(p)$.

	\item Fixing~$C$, suppose that $G\subseteq \QQ^\task(C)$ is $C$-generic; let $W = W^\task(C)[G]$, and let $\bar F$ be the generic $\mu$-sequence of sets. 

	 For a sequence $\bar x = \seq{x_\alpha}_{\alpha\in W}$ with $x_\alpha\in F_\alpha$, $\SS^\task(C,G,\bar x)$ is the notion of forcing consisting of conditions $(d,y)$, where
	\begin{itemize}
		\item $d,y\in 2^{<\lambda}$ and $\lth{d}= \lth{y}+1$; 
		\item $d$ is a closed subset of~$\lth{d}$; 
		\item for all $\alpha \in d \cap W$, $y\cap \alpha\notin F_\alpha\setminus \{{x_\alpha}\}	$;
		\item extension is extension in both coordinates. 
	\end{itemize}
\end{orderedlist}

Observe that indeed $\s^\task(C)$ is  $\lambda$-sparse: for a canonical sparse upper bound of an increasing sequence $\bar p = \seq{p_i}$ we choose a condition~$p^*$ defined by $\s(p^*) = \s(\bar p)\conc 0$, and $\bar F(p^*) = \bar F(\bar p) = \bigcup_i \bar F(p_i)$. Of course one of the main points is that $\alpha = |\s(\bar p)|$ is excluded from $\s(p^*)$, so we do not need to define $F_\alpha(p^*)$. Note that in fact, $\QQ^\task(C)$ is $\lllambda$-closed. Also observe that $\SS^\task(C,G,\bar x)$ is explicitly closed outside~$W$, by taking $\delta(d,y) = \lth{y}$. 

\begin{remark} \label{rmk:what_if_the_oracle_lies}
	In the example above, what is $\SS^{\task}(C,G,A)$ if~$A$ is not of the form~$\bar x$ where $x_{\alpha}\in F_{\alpha}(C)$? Note that the collection of ``legal'' $A$ is $\Pi^0_1(\HH_\lambda;C,G)$, so the generic definition $\SS^\task(C,G,A)$ can involve a check and output a trivial notion of forcing when~$A$ does not have the right form.
\end{remark}

\subsection{The task iteration} 
\label{sub:the_task_iteration}

We now define the \emph{$\lambda$-task iteration}. This will be a $\lambda$-sparse iteration $\seq{\PP_{\zeta},\QQ_{\zeta},\s_\zeta}$ of length $\lambda^+$. We let $\seq{\tau_\zeta}$ be the associated sequence, and as above use the abbreviations $W_\zeta = W_{\s_\zeta}$ and $U_\zeta = W_{\tau_\zeta}$. Together with the iteration we will define bookkeeping objects $\II,\ff$ such that $\II\subseteq \lambda^+$, and for all $\zeta<\lambda^+$,
\begin{orderedlist}
	\item If $\zeta\in \II$, then $\ff(\zeta)\in V(\PP_\zeta)$ is a $\lambda$-task~$\task$, and $\QQ_\zeta = \QQ^\task(\lambda\setminus U_\zeta)$, $\s_\zeta = \s^\task(\lambda\setminus U_\zeta)$; 

	\item If $\zeta\notin \II$, then $\ff(\zeta) = (\upsilon,A)$ where:
	\begin{itemize}
		\item $\upsilon<\zeta$ and $\upsilon\in \II$, and
		\item $A\in \PPP(\lambda)(\PP_\zeta)$;
	\end{itemize}
	and $\QQ_\zeta = \SS^\task(\lambda\setminus U_\upsilon,G_\upsilon,A)$, where $\task = \ff(\upsilon)$ and $G_\upsilon\in V(\PP_{\upsilon+1})$ is the generic for $\QQ_\upsilon = \QQ^\task(\lambda\setminus U_\upsilon)$. In this case, as $\QQ_\zeta$ is explicitly closed outside $W_\upsilon$, we let $\delta_\zeta$ witness this fact, and let $\s_\zeta$ be the associated explicit name for the empty set (see \cref{lem:closed_implies_sparsely_closed}). Note that $W_\upsilon \subseteq U_\zeta$ modulo clubs, so $\QQ_\zeta$ is explictly closed outside $U_\zeta$ as required. 
\end{orderedlist}

\begin{lemma} \label{lem:first_iteration:chain_condition}
	$\PP_{\lambda^+}$ is strategically $\lambda$-closed, and has the $\lambda^+$-chain condition. 
\end{lemma}

\begin{proof}
	Strategic $\lambda$-closure follows from \cref{lem:preservation_of_strategic_S_closure}. For the $\lambda^+$-chain condition, note that for all $\zeta<\lambda^+$, $\force_{\PP_\zeta} |\QQ_\zeta|\le \lambda$, and that the iteration is $\lllambda$-support. 
\end{proof}

\begin{corollary} \label{cor:first_iteration:GCH_and_preservation}
	$\PPP(\lambda)(\PP_{\lambda^+}) = \bigcup_{\zeta<\lambda^+} \PPP(\lambda)(\PP_\zeta)$; $\PP_{\lambda^+}$ is $\lllambda$-distributive; in $V(\PP_{\lambda^+})$, no cofinalities are changed and no cardinals are collapsed; if  $\GCH$ holds in~$V$ then it also holds in~$V(\PP_{\lambda^+})$.
\end{corollary}

\Cref{cor:first_iteration:GCH_and_preservation} allows us to define the bookkeeping functions so that:
\begin{itemize}
	\item For every $\lambda$-task $\task\in V(\PP_{\lambda^+})$ there is some $\xi\in \II$ such that $\ff(\xi) = \task$; 
	\item For every $\xi\in \II$, for every $A\in \PPP(\lambda)(\PP_{\lambda^+})$, there are unboundedly many $\zeta>\xi$ such that $\ff(\zeta) = (\xi,A)$. 
\end{itemize}

\subsection{Exact diamonds in the extension} 

\begin{proposition} \label{prop:first_iteration:diamonds_in_extension}
	In $V(\PP_{\lambda^+})$, for every regular $\theta<\lambda$, for every cardinal $\mu\in [2,\lambda)$, there is a stationary set $W\subseteq \Cof_\lambda(\theta)$ and a $\mu$-exact diamond sequence on~$W$. 
\end{proposition}

To prove this proposition, we let $\task = \task_{\Diamond}(\mu,\theta)$ be the obvious modification of the task $\task_\Diamond(\mu)$ above, where the requirement $\s(p)\subseteq C$ is replaced by $\s(p)\subseteq C\cap \Cof(\theta)$. 

We fix some $\upsilon\in \II$ such that $\ff(\upsilon) = \task$. Let $W = W_\upsilon = W^\task(\lambda\setminus U_\upsilon)[G_\upsilon]$ and let $\bar F = \bar F[G_\upsilon]$. 

The easier part is to check exactness.

\begin{lemma} \label{lem:first_iteration:diamonds:exact}
	In $V(\PP_{\lambda^+})$, for every $\bar x= \seq{x_\alpha}_{\alpha\in W}$ for which $x_\alpha\in F_\alpha$ for all $\alpha\in W$, 
	there are $\lambda^+$-many $Y\in \PPP(\lambda)$ which are not guessed by $\seq{F_\alpha\setminus \{x_\alpha\}}$. 
\end{lemma}

\begin{proof}
	Fix a set $T\subset \PPP(\lambda)(\PP_{\lambda^+})$ with $|T|\le \lambda$. Find some $\zeta$ sufficiently large so that $T\in V(\PP_\zeta)$ and $\ff(\zeta) = (\upsilon,\bar x)$, so $\QQ_\zeta = \SS^\task(\lambda\setminus U_\upsilon, G_\upsilon,\bar x)$. This notion of forcing adds a set~$Y\in \PPP(\lambda)\setminus T$ and a club~$D$ witnessing that $\seq{F_\alpha\setminus \{x_\alpha\}}$ does not guess~$Y$. 
\end{proof}

What is trickier is showing that $\bar F$ is a diamond sequence in $V(\PP_{\lambda^+})$. Certainly it is a diamond sequence in $V(\PP_{\upsilon+1})$: 
But we need to show that it remains a diamond sequence in $V(\PP_{\lambda^+})$. The argument follows, to a large degree, the main argument of~\cite{Sh:122}. 

\smallskip

Let $Z\in 2^\lambda(\PP_{\lambda^+})$, and let $C\subseteq \lambda$ in $V(\PP_{\lambda^+})$ be a club; find some limit $\xi\in (\upsilon,\lambda^+)$ such that $C,Z\in V(\PP_\xi)$. 

For the rest of the proof, we work in $V(\PP_\upsilon)$. In that universe, $U_\upsilon$ is co-fat, and $\PP_\xi/\PP_\upsilon$ is the result of the named $\lambda$-iteration $(\PP_\zeta/\PP_\upsilon,\QQ_\zeta,\s_\zeta)_{\zeta\in [\upsilon,\xi)}$ (\cref{lem:the_iteration_lemma}). 

Let 
 \[
 Q = \left\{ \zeta\in (\upsilon,\xi) \,:\, \ff(\zeta) = (\upsilon,\bar x) \text{ where in  }V(\PP_\zeta), \,\bar x\in \prod_{\alpha \in W} F_\alpha \right\}.
 \]
 For $\zeta\in Q$ we write $\bar x^\zeta = \seq{x^\zeta_\alpha}_{\alpha\in W}$ where $\ff(\zeta)= (\upsilon,\bar x^\zeta)$, and we let $(D_\zeta,Y_\zeta)\in V(\PP_{\zeta+1})$ be the pair of club and subset of~$\lambda$ which are added by~$\QQ_\zeta = \SS^\task(\lambda\setminus U_\upsilon,G_\upsilon,\bar x^\zeta)$. 

We note:
\begin{itemize}
	\item If $\zeta \in (\upsilon,\xi)\setminus Q$, then $\s_\zeta$ is $W$-sparse. 
\end{itemize}
For if $\zeta\in \II$ then $\s_\zeta$ is $\lambda$-sparse; if $\zeta \notin \II$ then $\QQ_\zeta$ is explicitly closed outside $W_\vartheta$ for some $\vartheta<\zeta$ distinct from~$\upsilon$; and $W_\vartheta\cap W = \emptyset$ modulo the club filter. So $\QQ_\zeta$ is explicitly $W$-closed; and $\s_\zeta$ (a $\QQ_\zeta$-name for the empty set) is $W$-sparse.

We apply \cref{cor:iteration_up_to_lambda_plus:existence_of_associated_sequence} to the iteration $\seq{\PP_\zeta/\PP_\upsilon}_{\zeta\in [\upsilon,\xi]}$, with $S = \lambda\setminus U_\upsilon$, obtaining an associated sequence of names $\smallseq{\tau_\zeta}_{\zeta\in (\upsilon,\xi]}$, which is sparse outside~$U_\upsilon$. However we modify the construction at successor steps $\zeta+1$ for $\zeta\in Q$ to determine not only the height $\delta_\zeta(p(\zeta)) = \smalllth{y^{p(\zeta)}}$ of a condition, but the actual values $(d,y)$ of the condition. So for all $p\in \dom \tau_\xi$, for all $\zeta\in \supp(p)\cap Q$, there is a pair $(d^p_\zeta, y^p_\zeta)\in V(\PP_\upsilon)$ (rather than $V(\PP_\zeta)$) such that $p\rest{\zeta}$ forces that $p(\zeta) = (d^p_\zeta, y^p_\zeta)$. Note that for such~$p$ and~$\zeta$ we have $\smalllth{d^p_\zeta} = \smalllth{\tau_\xi(p)}$. For brevity, let $\tau = \tau_\xi$.

\smallskip

There are two cases.

\smallskip

In the first case, we work above a condition $p^*\in \PP_\xi/\PP_\upsilon$ which forces that for all $\zeta\in Q$, $Z\ne Y_\zeta$.  Because $\PP_\xi/\PP_\upsilon$ is $\lllambda$-distributive, every $p\ge p^*$ in $\PP_\xi/\PP_\upsilon$ has an extension $r\in \dom \tau$ such that for some string $\pi\in 2^{<\lambda}$ with $\lth{\pi} > \smalllth{\tau(p)}$, 
\begin{orderedlist}
	\item $r\force \pi \prec Z$; and
	\item for every $\zeta\in Q\cap \supp(p)$, $\pi\perp y^r_\zeta$. 
\end{orderedlist}

Since $U_\upsilon$ is co-fat, we use our usual technique to get a sparse sequence $\bar p = \seq{p_i}_{i <  \theta}$\footnote{Recall that~$\theta$ is the cofinality on which $W$ concentrates.} in $\dom \tau$ such that for limit $i\le \theta$, $\sup_{j<i} \gamma_j\notin U_\upsilon$; and such that each $p_{i+1}$ separates~$Z$ from $y^{p_i}_\zeta$ as above, namely: for some $\pi_i$ of length in $[\gamma_i,\gamma_{i+1})$, we have $p_{i+1}\force \pi_i\prec Z$, and for all $\zeta\in Q\cap \supp(p_i)$, $\pi_i\perp y^{p_{i+1}}_\zeta$. Naturally $\seq{\pi_i}$ is increasing. Also, we ensure that $p_{i+1}$ forces that $C\cap [\gamma_i,\gamma_{i+1})$ is nonempty. We may assume that $\upsilon\in \supp(p_0)$. 

We are after a condition~$q$ which forces that~$Z$ is guessed by~$\bar F$ at some point in~$C$. Informally, we define~$q$ to be a modification of the canonical sparse upper bound~$p_\theta$ of~$\bar p$; the difference is in $q(\upsilon)$, which recall is a condition in $\QQ_\upsilon = \QQ^\task(\lambda\setminus U_\upsilon)$. Let $\alpha^* = \sup_{i<\theta} \gamma_i$, so $\cf(\alpha^*)=\theta$ and $\alpha^*\notin U_\upsilon$. 
Formally, by induction on~$\zeta\in (\upsilon,\xi]$, we define $q\rest{\zeta}$, extending each~$p_i\rest\zeta$, as follows. We declare that $\supp(q) = \bigcup_{i<\theta} \supp(p_i)$. For $\zeta\in \supp(q)\cap Q$, let $y^*_\zeta = \bigcup_{i<\theta} y^{p_i}_\zeta$ and $d^*_\zeta = \bigcup_{i<\theta} d^{p_i}_\zeta$. Also let $\pi^* = \bigcup_{i<\theta} \pi_i$.
\begin{equivalent}
	\item First we define $q(\upsilon)$ by letting $\s(q(\upsilon)) = \bigcup_{i<\theta}\s(p_i(\upsilon))\conc 1$ and $\bar F(q(\upsilon)) = \bigcup_{i<\theta}\bar F(p_i(\upsilon))\conc F_{\alpha^*}$, where $\pi^*\in F_{\alpha^*}$ but for all $\zeta\in \supp(q)\cap Q$, $y^*_\zeta\notin F_{\alpha^*}$. Note that this is a legitimate condition because $\alpha^*\in \Cof_\lambda(\theta)\setminus U_\upsilon$. 
\end{equivalent}
Let $\zeta>\upsilon$ be in $\supp(q)$, and suppose that $q\rest{[\upsilon,\zeta)}$ has already been defined, and that for all $i<\theta$, $q\rest{[\upsilon,\zeta)}$ extends $p_i\rest{[\upsilon,\zeta)}$. 
\begin{equivalent}[resume]
	\item Suppose that $\zeta\notin Q$. Since $q\rest{[\upsilon,\zeta)}$ extends each $p_i\rest{[\upsilon,\zeta)}$, it forces that a tail of $\seq{p_i(\zeta)}$ is a sparse sequence for $\s_\zeta$ (we use \ref{item:lambda_plus_iteration:taking_canonical_bounds_at_every_location} of \cref{def:associated_tau_sequence_regularity_properties:up_to_lambda_plus}). Also, $q\rest{[\upsilon,\zeta)}$ forces that $\alpha^*\in W$, so it forces that $\seq{p_i(\zeta)}$ has an upper bound, which we set to be $q(\zeta)$. 

	\item Suppose that $\zeta\in Q$. We declare that $q(\zeta) = (d^*_\zeta\conc 1,y^*_\zeta)$. This is a legitimate condition because $\alpha^*\in W$ and $y_\zeta^*\notin F_{\alpha^*}$. 
\end{equivalent}

The condition~$q$ forces that $\pi^* = Z\cap \alpha^*$, $\alpha^*\in C\cap W$ and $\pi^*\in F_{\alpha^*}$, finishing the proof in this case.

 \medskip

In the second case, we work above a condition $p^*$ which forces that $Z = Y_{\varrho}$ for some $\varrho\in Q$. (Note that we still have to work in $V(\PP_\xi)$ rather than $V(\PP_{\varrho})$, as $C\in V(\PP_\xi)$ may not be in $V(\PP_{\varrho})$.) We may assume that $p^*\in \dom \tau$ and $\varrho\in\supp(p^*)$. Let $\alpha = \lth{\tau(p^*)}+\mu$, where, recall, $\mu$ is the size of each $F_\beta$. Let $p^*(\varrho) = (d^*,y^*)$. Find strings $y^i$, for $i<\mu$, of length~$\alpha$, each extending $y^*$, which are pairwise incomparable; let $p^i$ agree with $p^*$ except that $p^i(\varrho) = (d,y^i)$, where $\lth{d}= \alpha+1$ is an extension of~$d^*$ by a string of zeros. 

Using \cref{lem:extending_compatibles_and_incompatibles}, we follow our usual construction, this time constructing an array $\seq{p^i_j}_{i<\mu,j< \theta}$ such that:
\begin{orderedlist}
	\item For each~$i<\mu$, $\seq{p^i_j}_{j\le \theta}$ is a sparse sequence for~$\tau$;
	\item For each $i<i'\le \theta$, $p^i_j\rest{\varrho} = p^{i'}_j\rest{\varrho}$, call the common value $r_j$;
	\item $p^i_0 = p^i$; 
	\item For each $i<\mu$ and $j< \theta$, $p^i_j\in A^\xi_{\varrho}$; 
	\item For each $i,i'<\mu$, $j<\theta$, and $\zeta\in Q\cap \supp(p^i_j)$ other than $\varrho$, $y^{p^i_{j+1}}_\zeta\perp y^{p^{i'}_{j+1}}_\varrho$. 
\end{orderedlist}
At step~$\theta$, instead of taking the sparse upper bound, we inductively define a condition~$q\in \dom \tau_\varrho$ by induction as above. For $i<\mu$ let $y^i_\varrho = \bigcup_{j<\theta} y^{p^i_j}_\varrho$. Let $\alpha^* = \lth{\tau_{\varrho}(\bar r)}$, which has cofinality~$\theta$ and is out of~$U_\upsilon$. We start with setting $q(\upsilon)$ so that $\alpha^*\in \s_\upsilon(q(\upsilon))$ and $F_{\alpha^*}$, according to $q(\upsilon)$, is $\left\{ y^i_\varrho  \,:\, i<\mu   \right\}$. We then define $q(\zeta)$ for $\zeta\in \supp(q) =  \bigcup_j \supp(r_j)$ as above, noting that for all $\zeta\in Q\cap \supp(q)$, for all $i<\mu$, $y^i_\varrho\ne y^q_\zeta$, where of course $y^q_\zeta = \bigcup_{j<\theta} y^{r_j}_\zeta$. 

Now, we extend~$q$ to a condition $q'\in \dom \tau_\varrho$ such that for some $\pi$ of length~$\alpha^*$, $q'\force x^\varrho_{\alpha^*}= \pi$. Since $q'$ decides the value of $F_{\alpha^*}$, necessarily $\pi=  y^{i^*}_\varrho$ for some $i^*<\mu$. As above, extend~$q'$ to a condition $s\in \dom \tau_\xi$ extending $p^{i^*}_j$ for all $j<\theta$. Defining $s(\varrho) = (d^{i^*}_\varrho,y^{i^*}_{\varrho})$ (where of course $d^i_\varrho = (\bigcup_{j<\theta}d^{p^i_j}_\varrho)\conc 1$) is legitimate since~$q'$ forces that $y^{i^*}_{\varrho}\notin F_{\alpha^*}\setminus \{x^\varrho_{\alpha^*}\}$. The condition~$s$ forces that $\alpha^*\in C$ and that $Z\cap{\alpha^*} = y^{i^*}_{\varrho}\in F_{\alpha^*}$.

\section{The task axiom} 
\label{sec:the_task_axiom}

\subsection{Named iterations beyond $\lambda^+$} 
\label{sub:named_iterations_beyond_lambda_plus}


Let $\UU$ be a family of subsets of~$\lambda$. We define the notion of forcing $\RR(\UU)$:
\begin{itemize}
	\item Conditions are pairs $p = (u(p),\s(p))$ where $\s(p)\in 2^{<\lambda}$ and $u(p)\subseteq \UU$ has size $<\!\lambda$.
	\item $q$ extends~$p$ if $u(p)\subseteq u(q)$, $\s(p)\preceq \s(q)$ and for all $U\in u(p)$, $U\cap [\lth{\s(p)},\lth{\s(q)})\subseteq \s(q)$. 
\end{itemize}

This notion of forcing is $\lllambda$-closed; $\s$ is an explicit $\RR(\UU)$-name for a subset of~$\lambda$. In $V(\RR(\UU))$, $W_\s$ is an upper bound of $\UU$ modulo clubs (indeed modulo bounded subsets of~$\lambda$). We will not need to directly appeal to the following lemma, since in our context we will use the heavier machinery of sparse names. 

\begin{lemma} \label{lem:R_U_preserves_stationarity_outside_U}
	For every $S\subseteq \lambda$ in~$V$, if for all $U\in \UU$, $S\setminus U$ is stationary, then in $V(\RR(\UU))$, $S\setminus W_\s$ is stationary. \qed
\end{lemma}

We can now give a general definition of named $\lambda$-iterations (of any length):

\begin{definition} \label{def:named_lambda_iteration:general}
	A sequence $\seq{\PP_\zeta}_{\zeta<\xi}$ is a named $\lambda$-iteration if it satisfies conditions (i)-(iii) of \cref{def:named_iterations}, together with:
	\begin{orderedlist}[start=4]
		\item For $\zeta<\xi$ such that $\cf(\zeta)>\lambda$, $\QQ_\zeta = \RR(\left\{ W_{\s_\upsilon} \,:\,  \upsilon<\zeta \right\})$ and~$\s_\zeta$ is the corresponding explicit name. 
	\end{orderedlist}
\end{definition}

Using this definition, as above, for such an iteration we can define a sequence $\seq{\uu_\zeta}_{\zeta<\xi, \cf(\zeta)\le \lambda}$ with $\uu_\zeta\in \PPP(\lambda)/\NS_\lambda(\PP_{\zeta})$. This will be an increasing sequence in $\PPP(\lambda)/\NS_\lambda$. If $\zeta = \upsilon+1$ and $\cf(\upsilon)>\lambda$ then we let $\uu_\zeta = [W_{\s_\upsilon}]_{\NS_\lambda}$. Otherwise we can take least upper bounds as above. \Cref{def:sparse_iteration} now makes sense for named $\lambda$-iterations of length beyond $\lambda^+$ as well. Similarly we can generalise \cref{def:associated_tau_sequence_regularity_properties:up_to_lambda_plus}; the only change is that $\tau_\zeta$ is defined only for $\zeta\in [1,\xi)\cap \Cof(\le \lambda)$. We can then extend \cref{prop:iteration_up_to_lambda_plus:existence_of_associated_sequence:inductive_step}:

\begin{proposition} \label{prop:obtaining_tau_for_long_iteration}
	Let $S\subseteq \lambda$ be fat. Suppose that $\seq{\PP_\zeta}_{\zeta\le \xi}$ is an $S$-sparse iteration, and that $\bar \tau = \seq{\tau_\zeta}_{\zeta\in [1,\xi)\cap \Cof(\le\!\lambda)}$ is an associated sequence for $\bar \PP\rest{\xi}$. Suppose that $\cf(\xi)\le \lambda$.  

	Then there is a name $\tau_\xi$ such that $\bar \tau\conc \tau_\xi$ is an associated sequence for $\bar \PP$. 
\end{proposition}

\begin{proof}
	The proof of \cref{prop:iteration_up_to_lambda_plus:existence_of_associated_sequence:inductive_step} holds in all cases handled in that proposition, so we assume that $\xi = \vartheta+1$ where $\cf(\vartheta)>\lambda$. In this case, by \cref{lem:preservation_of_strategic_S_closure}, $\PP_{\vartheta}$ is $\lllambda$-distributive.

	For brevity, let $I = \xi\cap \Cof(\le \lambda)$. Also for brevity, let $\RR = \QQ_{\vartheta} = \RR(\left\{ W_{\s_\zeta} \,:\,  \zeta<\xi \right\})$. Since $\PP_{\vartheta}$ is $\lllambda$-distributive, the set of conditions 
	$(p,q)\in \PP_\vartheta\starr \RR = \PP_{\xi}$ such that for some $v\subseteq I$ of size $<\!\lambda$ and some string $\pi\in 2^{<\lambda}$, 
	\[
		p\force_{\PP_\vartheta}  \s(q) = \pi \andd u(q) = \left\{  W_{\s_\zeta}  \,:\, \zeta\in v   \right\}
	\]
	is dense in $\PP_\vartheta\starr \RR$. We restrict ourselves to this dense subset, and write $v(q), \s(q)$ for the corresponding $v$ and~$\pi$. 

	\smallskip
	
	Let $\langle{A^\beta_\alpha}\rangle$ witness the weak coherence of $\bar \tau$. We let $(p,q)\in \dom \tau_{\xi}$ if:
		\begin{orderedlist}
		\item \label{item:dom_lambda_plus:where_p_lives}
		 $p\in \PP_{\sup v(q)}$;
		 \item \label{item:dom_lambda_plus:add_omega}
		 For all $\delta$, $\delta\in v(q) \Iff \delta+1\in v(q)$; 
		\item \label{item:dom_lambda_plus:restriction_in_dom_tau}
		for all $\delta\in v(q)$, $p\rest{\delta}\in \dom \tau_\delta$ and $\lth{\tau_\delta(p\rest{\delta})} = \lth{\s(q)}$;  and
		\item \label{item:dom_lambda_plus:coherence}
		for all $\delta\in v(q)$, for all $\epsilon<\delta$ in~$I$, $p\rest{\delta}\in A^\delta_\epsilon$ if and only if $\epsilon\in v(q)$. 
	\end{orderedlist}
	For $(p,q)\in \dom \tau_{\xi}$ we let $\tau_{\xi}(p,q) = \s(q)$. If $(\bar p,\bar q)= \seq{(p_i,q_i)}_{i<i^*}$ is an increasing sequence from $\dom \tau_\xi$, then we let $\cb_{\tau_\xi}(\bar p,\bar q) = (p^*,q^*)$ if:
	\begin{orderedlist}[resume]
		\item $(p^*,q^*)\in \dom \tau_\xi$;
		\item $v(q^*) = \bigcup_{i<i^*} v(q_i)$; 
		\item $\s(q^*) = \bigcup_{i<i^*}\s(q_i)\conc 0$; 
		\item For all $\delta\in v(q^*)$, $p^*\rest{\delta}$ is $\cb_{\tau_\delta}$ of a tail of $\bar p\rest{\delta}$. 
	\end{orderedlist}

	\smallskip
	
	We show that $\dom \tau_\xi$ is dense in $\PP_\xi$. Given $(p_0,q_0)\in \PP_\xi$, we obtain a $\lambda$-filtration~$\bar N$ and a sequence $\seq{\gamma_n}_{n<\w}$ with limit $\gamma_\w\in S$ as usual. We define an increasing sequence $\seq{p_n,q_n}_{n<\w}$ with $(p_n,q_n)\in N_{\gamma_n}$. We require that:
	\begin{itemize}
		\item $\lth{\s(q_n)}> \gamma_{n-1}$;
		\item $N_{\gamma_{n-1}}\cap I \subseteq v(q_n)$; 
		\item $p_n\in \PP_{\zeta_n}$ for some $\zeta_n\in N_{\gamma_n}$, $\zeta_n> \sup v(q_n)$;
		\item $\lth{\tau_{\zeta_n}(p_n)}> \gamma_{n-1}$; and
		\item for all $\delta\in v(q_n)$, $p_n\in A^{\zeta_n}_{\delta}$. 
	\end{itemize}
	Note that having size $<\!\lambda$,  $v(q_n)\subset N_{\gamma_n}$.

	Let $v^* = \bigcup_n v(q_n) = N_{\gamma_\w}\cap I$, which is cofinal in $N_{\gamma_\w}\cap \xi$; let $\zeta^* = \sup v^*$. Also let $\pi^* = \bigcup_n \s(q_n)$, which has length $\gamma_\w$.  Since $|v^*| < \lambda$, $\PP_{\zeta^*}$ is the inverse limit of $\seq{\PP_\delta}_{\delta\in v^*}$. For $\delta\in v^*$, a tail of $\seq{p_n\rest{\delta}}$ is in $\dom \tau_\delta$; we let $r_\delta = \cb_{\tau_\delta}(\seq{p_n\rest{\delta}})$ be its canonical upper bound. By coherence, $r_\delta\rest{\epsilon} = r_\epsilon$ for $\epsilon<\delta$ in $v^*$, 
	so we obtain an inverse limit $r\in \PP_{\zeta^*}$; this condition extends each~$p_n$. Define $(r,s)\in \PP_{\vartheta}\starr \RR$ by letting $\s(s) = \pi^*\conc 1$, and $v(s)= v^*$. We claim that $(r,s)\in \dom \tau_\xi$. Requirements \ref{item:dom_lambda_plus:where_p_lives} -- \ref{item:dom_lambda_plus:restriction_in_dom_tau} and one direction of \ref{item:dom_lambda_plus:coherence} are clear; for the other direction of \ref{item:dom_lambda_plus:coherence}, we suppose that $\delta\in v^*$, $\epsilon<\delta$, and $r_\delta = r\rest \delta\in A^\delta_\epsilon$; we need to show that $\epsilon\in v^*$. By the continuity property~\ref{item:weak_coherence:continuity} of \cref{def:weak_coherent_sequence}, for large enough~$n$, $p_n\rest{\delta}\in A^\delta_\epsilon$. Fixing some large~$n$, we know (by \ref{item:weak_coherence:number_bound}) that~$\epsilon$ is among the fewer than~$\lambda$-many $\epsilon'<\delta$ such that $p_n\rest{\delta}\in A^{\delta}_{\epsilon'}$. Since $\seq{A^\beta_\alpha},I,p_n\in N_{\gamma_n}$ it follows that each such~$\epsilon'$ is in $N_{\gamma_n}$, so $\epsilon\in I\cap N_{\gamma_n}\subseteq v(q_{n+1})\subset v^*$. 

	To show that $(r,s)$ extends $(p_n,q_n)$, we need to show that if $\delta\in v(q_n)$ and $\alpha\in (\lth{\s(q_n)},\gamma_\w]$, then $r$ forces that if $\alpha\in W_{\s_\delta}$, then $\s(s)(\alpha)=1$. But we may assume that $\alpha<\gamma_\w$; if $m>n$ is sufficiently large so that $\lth{\s(q_m)}, \lth{\tau_\delta(p_m\rest{\delta})}>\alpha$, then $p_m$ forces what's required, as it forces that~$q_m$ extends~$q_n$. 

	\medskip
	
	We check that the requirements of \cref{def:associated_tau_sequence_regularity_properties:up_to_lambda_plus} hold. We first need to check that $\tau_\xi$ is $S$-sparse. Suppose that $(\bar p, \bar q) = \seq{p_i,q_i}_{i<i^*}$ is a sparse sequence for~$\tau_\zeta$, with $\alpha^* = \lth{\tau_\zeta(\bar p,\bar q)}\in S$. Let $v^* = \bigcup_{i<i^*} v(q_i)$; let $\pi^* = \bigcup_{i<i^*} \s(q_i)$. Let $\zeta^* = \sup v^*$. For all $\delta\in v^*$, a tail of $\bar p\rest{\delta}$ is in $\dom \tau_\delta$ and is a sparse sequence for $\tau_\delta$; we let $r_\delta$ be the canonical upper bound; we let~$r$ be the inverse limit of $\seq{r_\delta}$; we define~$s$ by $\s(s)= \pi^*\conc 0$ and $v(s) = v^*$. It is not difficult to see that $(r,s)\in \dom \tau_\xi$ and $(r,s) = \cb_{\tau_\xi}(\bar p,\bar q)$, except that we need to show that $(r,s)$ extends each $(p_i,q_i)$. The only potentially contentious point is the legitimacy of setting $\s(s)(\alpha^*)=0$. Let $\delta\in v(q_i)$. Then $\delta+1\in v(q_i)\subseteq v(s)$, whence $p_i\rest{\delta+1}\in \dom \tau_{\delta+1}$. By \ref{item:lambda_plus_iteration:uniform_length_of_conditions} for $\tau_{\delta+1}$, $\delta\in \supp(p_i)$, and since $\tau_{\delta+1}(r_{\delta+1})(\alpha^*)=0$, by \ref{item:lambda_plus_iteration:subset_relations_from_a_point} (applied to $q = r_{\delta+1}$ and $p = p_i\rest{\delta+1}$), $\s_\delta(r)(\alpha^*)=0$, so $\delta$ and~$\alpha^*$ are not an obstacle to $(r,s)$ extending $(p_i,q_i)$.

	\smallskip
	
	\ref{item:lambda_plus_iteration:U_s_correct} follows from the fact that we defined $\uu_\xi = [W_{\s_\vartheta}]_{\NS_\lambda}$.

	\smallskip

	For \ref{item:lambda_plus_iteration:weak_coherence}, for $\delta\in I$ and $(p,q)\in \dom \tau_\xi$, we let $(p,q)\in A^\xi_\delta$ if and only if $\delta\in v(q)$. The requirements of \cref{def:weak_coherent_sequence} follow from our definitions. 

	For $\ref{item:lambda_plus_iteration:uniform_length_of_conditions}$, let $\zeta\in \supp(p)$; since $p\in \PP_{\sup v(q)}$, $\zeta<\sup v(q)$; let $\delta\in v(q)$, $\delta>\zeta$. Then $p\rest{\delta}\in \dom \tau_\delta$, and so $p\rest{\delta}\in A^\delta_{\zeta+1}$; it follows that $\zeta+1\in v(q)$. On the other hand we note that $\vartheta\in \supp (p,q)$. 

	\ref{item:lambda_plus_iteration:zero_in_sigmas_implies_in_tau} is immediate: if for all $\delta\in \supp(p,q)$, $\s_\delta(p,q)(\alpha)=0$, then $\tau_\xi(p,q)(\alpha)= \s_{\vartheta}(q)(\alpha)=0$. 

	For \ref{item:lambda_plus_iteration:subset_relations_from_a_point}, suppose that $\zeta\in \supp(p,q)$, and that $(p',q')$ extends $(p,q)$. We may assume that $\zeta<\vartheta$, so $\zeta\in v(q)$. Then $\s_\zeta(p'(\zeta))\setminus \lth{\s_\zeta(p(\zeta))}\subseteq \s_\vartheta(q') = \tau(p',q')$ because~$p'$ forces that~$q'$ extends~$q$. 

	Finally, \ref{item:lambda_plus_iteration:continuity_of_support_in_canonical_bounds} and \ref{item:lambda_plus_iteration:taking_canonical_bounds_at_every_location} follow from our definitions. 
\end{proof}

\subsubsection{Quotients of named iterations} 

Let~$S$ be fat, and let $\seq{\PP_{\zeta}}_{\zeta<\xi}$ be an $S$-sparse iteration. Let $\upsilon<\xi$. If $\cf(\upsilon)\le \lambda$ then the usual argument shows that the quotient iteration $\seq{\PP_\zeta/\PP_\upsilon}_{\zeta\in [\upsilon,\xi)}$ is, in $V(\PP_\upsilon)$, an $S\setminus U_\upsilon$-sparse iteration. 

What happens though if $\cf(\upsilon)>\lambda$? In this case, the iteration $\seq{\PP_\zeta/\PP_\upsilon}$ is a one-step extension $\RR(\UU)\starr \seq{\PP_{\zeta}/\PP_{\upsilon+1}}_{\zeta \in (\upsilon,\xi)}$, where $\UU = \left\{ W_\zeta \,:\,  \zeta<\upsilon \right\}$. The second step, as mentioned, is an $S\setminus U_{\upsilon+1}$-sparse iteration, where $U_{\upsilon+1}= W_\upsilon$ comes from the generic for $\RR(\UU)$; we know that $S\setminus U_{\upsilon+1}$ is fat in $V(\PP_{\upsilon+1})$. 

\subsection{Nice tasks and the task axiom} 
\label{sub:nice_tasks}

\begin{definition} \label{def:a_task_respecting_iteration}
	Let $\task$ be a $\lambda$-task, let $S\subseteq \lambda$ be fat. An $S$-sparse iteration $\seq{\PP_\zeta}_{\zeta<\xi}$ 
	is \emph{$\task$-friendly} if $\QQ_0 = \QQ^\task(S)$, and letting $G_0$ be the generic for $\QQ_0$ and $W = W^\task(S)[G_0] = W_0$, for all $\zeta\in [1,\xi)$ with $\cf(\zeta)\le \lambda$, in $V(\PP_\zeta)$,
	\begin{orderedlist}
		\item either $\s_\zeta$ is $W$-sparse; or
		\item $\QQ_\zeta = \SS^\task(S, G_0,A)$ for some~$A$. 
	\end{orderedlist}
\end{definition}

For example, for the task iteration above, for every $\upsilon\in \II$, in $V(\PP_\upsilon)$, the quotient iteration $\seq{\PP_\zeta/\PP_{\upsilon}}$ is $\ff(\upsilon)$-friendly (with $S = \lambda\setminus U_\upsilon$).

\begin{definition} \label{def:correctness_condition}
	Let $\task$ be a $\lambda$-task. A \emph{correctness condition} for~$\task$ is a $\Pi^1_1$ sentence~$\psi$ such that for all fat sets~$C$, $\QQ^\task(C)$ forces that $(\HH_\lambda;O^\task,C,G)\models \psi$ (where~$O^\task$ is the oracle for~$\task$). 
\end{definition}

For example, for $\task = \task_\Diamond(\mu)$ (or $\task_\Diamond(\mu,\theta)$) as above, the correctness condition is that the resulting sequence $\seq{F_{\alpha}}_{\alpha\in W}$ is a $\mu$-diamond sequence.\footnote{Note that if~$G\subset \QQ^\task(C)$ is merely $C$-generic, then for this task it is not the case that $\psi$ holds in $(\HH_\lambda;C,G)$; for~$\psi$ to hold we need a filter~$G$ fully generic over~$V$.}

\begin{definition} \label{def:nice_task}
	Let~$\task$ be a $\lambda$-task and let~$\psi$ be a correctness condition for~$\task$. We say that the pair $(\task,\psi)$ is \emph{nice} if for any $\lllambda$-distributive~$\RR$, in $V(\RR)$, for any fat~$S$, for any $S$-sparse iteration $\seq{\PP_\zeta}_{\zeta\le \xi}$ with  $\xi <\lambda^{++}$ which is $\task$-friendly, in $V(\RR\starr \PP_\xi)$, $(\HH_\lambda;O,S,G_0)\models \psi$ (where as above $G_0$ is the generic for~$\QQ_0$). 
\end{definition}

\begin{definition} \label{def:satisfying_a_task}
	Let $\task$ be a $\lambda$-task and let~$\psi$ be a correctness condition for~$\task$. We say that the pair $(\task,\psi)$ is \emph{satisfied} if:
	\begin{equivalent}
		\item There is a fat set~$C$ and a filter $G\subset \QQ^\task(C)$ which is $(O^\task,C)$-generic; 
		\item $(\HH_\lambda;O^\task,C,G)\models \psi$; and 
		\item For all $A\in \PPP(\lambda)$, there is a filter $H_A\subset \SS^\task(C,G,A)$ which is $(O^\task,C,G,A)$-generic. 
	\end{equivalent}
\end{definition}

We can now define the $\lambda$-task axiom. 

\begin{definition} \label{def:lambda_task_axiom}
	The $\lambda$-task axiom $\TaskAxiom{\lambda}$ states:
	\begin{itemize}
		\item[] Every nice pair $(\task,\psi)$ is satisfied.
	\end{itemize}
\end{definition}

Again we emphasise that throughout, we assume that  $\GCH$ holds below~$\lambda$. 

\subsubsection{Diamonds and the task axiom} 

As mentioned above, for $\task = \task_\Diamond(\mu,\theta)$, the natural correctness condition~$\psi$ for~$\task$ is that $\bar F = \seq{F_\alpha}_{\alpha\in W}$ is a $\mu$-diamond sequence. Further, the argument proving \cref{prop:first_iteration:diamonds_in_extension} shows that $(\task,\psi)$ is nice. We modify the interpretation of the extra oracle~$A$ in the definition of $\SS^\task(C,G,A)$, so that $A$ not only codes the sequence $\bar x = \seq{x_\alpha}_{\alpha\in W}$ but also codes a family of $\lambda$-many $Z\in \PPP(\lambda)$, the result being that the generic~$Y$ satisfies $Y\ne Z$ for all~$Z$ coded by~$A$. The result is:

\begin{proposition} \label{prop:task_axiom_and_diamond}
	$\TaskAxiom{\lambda}$ implies that for all $\mu<\lambda$ and all regular $\theta<\lambda$, there is an exact $\mu$-diamond sequence concentrating on $\Cof_\lambda(\theta)$. 
\end{proposition}

\subsection{Consistency of the task axiom} 

We remind the reader that as usual, $\lambda$ is regular and uncountable, and that  $\GCH$ holds below~$\lambda$. In this subsection we prove:

\begin{proposition} \label{prop:task_axiom_consistency:single_cardinal}
	There is a notion of forcing $\PP$ which is $\lllambda$-distributive and has the $\lambda^+$-chain condition (and so preserves cardinals and cofinalities), preserves  $\GCH$ (if it holds in~$V$), and such that the task axiom $\TaskAxiom{\lambda}$ holds in $V(\PP)$. 
\end{proposition}


\begin{proof}
	Toward constructing~$\PP$, we define a $\lambda$-sparse iteration $\seq{\PP_{\zeta}}_{\zeta<\lambda^{++}}$; $\PP$ will be a proper initial segment of this iteration. As above we have book-keeping devices $\II$ and~$\ff$, except that for $\upsilon\in \II$, $\ff(\upsilon)$ is a pair $(\task,\psi)$ consisting of a task and a correctness condition for that task. Also, of course, $\ff$ is only defined on $\lambda^{++}\cap \Cof(\le \lambda)$, as $\QQ_\zeta$ is prescribed when $\cf(\zeta) = \lambda^+$. Also, we will repeat tasks when earlier attempts have resulted in failure. Suppose that $\ff(\upsilon) = (\task,\psi)$, and that $\zeta>\upsilon$. We say that \emph{$\ff(\upsilon)$ has failed by stage~$\zeta$} if there is some condition in $\PP_\zeta$ which forces that in $V(\PP_\zeta)$, $(\HH_\lambda;O,\lambda\setminus U_\upsilon,G_\upsilon)\models \lnot \psi$. 
We keep our books so that:
\begin{itemize}
	\item For every pair $(\task,\psi)\in V(\PP_{\lambda^{++}})$, either
	\begin{itemize}
		\item There is some $\upsilon\in \II$ such that $\ff(\upsilon) = (\task,\psi)$ and for all $\zeta>\upsilon$, $\ff(\upsilon)$ does not fail by stage~$\zeta$; or
		\item There are unboundedly many $\upsilon\in \II$ such that $\ff(\upsilon) = (\task,\psi)$, and if $\upsilon<\upsilon'$ are in~$\II$ and $\ff(\upsilon) = \ff(\upsilon')$, then $\ff(\upsilon)$ has failed by stage~$\upsilon'$. 
	\end{itemize}
	\item In the first case, with~$\upsilon$ last such that $\ff(\upsilon)= (\task,\psi)$, for all $A\in \PPP(\lambda)(\PP_{\lambda^{++}})$, there are unboundedly many $\zeta<\lambda^{++}$ such that $\ff(\zeta) = (\upsilon,A)$.
 \end{itemize} 

Then by a standard closure argument, we can find some $\xi^*<\lambda^{++}$ such that:
\begin{orderedlist}
	\item $\cf(\xi^*) = \lambda^+$;
	\item For every pair $(\task,\psi)\in V(\PP_{\xi^*})$, if there is a last $\upsilon\in \II$ such that $\ff(\upsilon)= (\task,\psi)$, then this $\upsilon$ is smaller than $\xi^*$;
	\item If~$\upsilon<\xi^*$ is last with $\ff(\upsilon) = (\task,\psi)$, then for all $A\in \PPP(\lambda)({\PP_{\xi^*}})$ there are unboundedly many $\zeta<\xi^*$ such that $\ff(\zeta) = (\upsilon,A)$.
\end{orderedlist}

We let $\PP = \PP_{\xi^*}$. Suppose that $\task\in V(\PP)$, that~$\psi$ is a correctness condition for~$\task$, and that $(\task,\psi)$ is nice in $V(\PP)$. Then there is a last~$\upsilon\in \II$ such that $\ff(\upsilon) = (\task,\psi)$. For otherwise, there is some $\upsilon>\xi^*$ such that $\ff(\upsilon) = (\task,\psi)$ and some $\xi>\upsilon$ such that $\ff(\upsilon)$ has failed by stage~$\xi$. Then the iteration $\seq{\PP_{\zeta}/\PP_{\upsilon}}_{\zeta\in [\upsilon^*,\xi)}$ witnesses (in $V(\PP_{\upsilon})$, which is an extension of $V(\PP_{\xi^*})$ by $\PP_\upsilon/\PP_{\xi^*}$, which 
as a quotient of a $\lllambda$-distributive notion of forcing is also $\lllambda$-distributive), that $(\task,\psi)$ is not nice.

Taking this last~$\upsilon$, we have $\upsilon<\xi^*$; we know that~$\psi$ holds in $V(\PP_{\xi^*})$; and our bookkeeping ensures that $(\task,\psi)$ is satisfied in~$V(\PP_{\xi})$. 
\end{proof}

Iterating with set Easton support, we get:

\begin{corollary} \label{cor:the_class_iteration}
	Assuming  $\GCH$, there is a class forcing extension preserving all cofinalities, cardinals, and $\GCH$, and in which $\TaskAxiom{\lambda}$ holds for every regular uncountable cardinal~$\lambda$. 
\end{corollary}

\subsection{Uniformisation} 
\label{sub:uniformisation}

Let $W\subseteq \lambda$. Recall that a \emph{ladder system} on~$W$ is a sequence $\bar E = \seq{E_\alpha}_{\alpha\in W}$ such that each $E_\alpha$ is an unbounded subset of~$\alpha$ of order-type $\cf(\alpha)$. A \emph{2-colouring} of a ladder system $\bar E$ is a sequence $\bar c = \seq{c_\alpha}_{\alpha\in W}$ such that for all $\alpha\in W$, $c_\alpha\colon E_\alpha \to 2$. A function $g\colon \lambda\to 2$ \emph{uniformises} the 2-colouring $\bar c$ if for every $\alpha\in W$, $g\rest{E_\alpha}=^* c_\alpha$, meaning that $\left\{ \gamma\in E_\alpha \,:\,  g(\gamma)\ne c_\alpha(\gamma) \right\}$ is bounded below~$\alpha$. 

A ladder system $\bar E$ on~$W$ \emph{has uniformisation} if every 2-colouring of $\bar e$ is uniformised by some~$g$.

Since we are assuming $\GCH$ below~$\lambda$, there are some restrictions on what kind of uniformisation we can deduce from the task axiom. If $\lambda = \mu^+$ and $\theta = \cf(\mu)$ then no ladder system on any stationary $W\subseteq \Cof_\lambda(\ne \theta)$ can have uniformisation, as $\Diamond_W$ holds \cite{Gregory,Sh:108,Sh:922}. Also, if $\lambda = \mu^+$ and~$\mu$ is singular, then for no stationary $W\subseteq \lambda$ is it the case that every ladder system on~$W$ has uniformisation \cite{Sh:667}.

\begin{proposition} \label{prop:uniformisation}
  Suppose that $\lambda$ is inaccessible and $\theta<\lambda$ is regular; or that $\lambda = \mu^+$ and $\theta = \cf(\mu)$. Then $\TaskAxiom{\lambda}$ implies that for any ladder system $\bar E$ on $\Cof_\lambda(\theta)$ there is some stationary $W\subseteq \Cof_\lambda(\theta)$ such that the restriction $\bar E\rest{W}$ has uniformisation. 
  \end{proposition}

\begin{proof}
Let $\bar E$ be a ladder system on $\Cof_\lambda(\theta)$. Define the following task $\task = \task_{\UP}(\theta)$:

\begin{orderedlist}
	\item For $C\subseteq \lambda$, $p\in \QQ^\task(C)$ if $p\in 2^{<\lambda}$, $p\subset C\cap \Cof(\theta)$, and the restriction of $\bar E$ to~$p$ has uniformisation. 

	We let $\sigma^\task(C)(p)  = p$. 

\item If $G\subseteq \QQ^\task(C)$ is $C$-generic, then we let $W = W[G] = W_{\s^\task(C)}[G]$. 

If  $\bar c$ is a 2-colouring of $\bar E\rest{W}$, then $\SS^\task(C,G,\bar c)$ consists of conditions $q\in 2^{<\lambda}$ which are an initial segment of a uniformising function: for all $\alpha\in W$ with $\alpha \le |q|$, $q\rest{E_\alpha} =^* c_\alpha$. 
\end{orderedlist}

First, we verify \cref{def:task}. It is clear that $\SS^\task(G,\bar c)$ is explicitly closed outside~$W$; take $\delta(q) = |q|$. To see that $\s^\task(C)$ is $\lambda$-sparse (using $\cb_\s(\bar p) = (\bigcup\bar p)\conc 0$), suppose that $\bar p = \seq{p_i}_{i<i^*}$ is sparse for~$\s$. Let $p^*=(\bigcup \bar p)\conc 0$. We need to argue that the restriction of~$\bar E$ to~$p^*$ has uniformisation. This follows from the closed set disjoint from~$p^*$ determined by the sequence $\bar p$: let $\gamma_i = |p_i|$. Every $\alpha\in p^*$ is in $(\gamma_i,\gamma_{i+1})$ for some~$i$, the salient point being that $\alpha > \gamma_i$. So if $\bar c$ is a 2-colouring of $\bar E\rest{p^*}$, then as each $p_i$ is in $\QQ^\task(C)$, we let $g_i\colon \gamma_i\to 2$ uniformise $\bar c\rest{p_i}$; define $g\colon |p^*|\to 2$ by letting~$g$ agree with~$g_{i+1}$ on ${[\gamma_i,\gamma_{i+1})}$.

\smallskip

Let~$\psi$ be the correctness condition which states that $W = W[G]$ is stationary. We argue that this is forced by $\QQ^\task(C)$. To see this, let $D\in V(\QQ^\task(C))$ be a club; and suppose that~$C$ is fat. Obtain a filtration $\bar N$ and a sequence $\seq{\gamma_i}_{i<\theta}$ (with limits in~$C$) as usual. Define an increasing sparse sequence $\bar p$ of length~$\theta$ as usual, so any upper bound forces that $\alpha^* = \sup_i \gamma_i$ is in~$D$; it is, of course, also in $C\cap \Cof(\theta)$. We argue that $p^* = (\bigcup \bar p)\conc 1$ is a valid condition in $\QQ^\task(C)$, namely, that it has uniformisation. Let $\bar c$ be a 2-colouring of $\bar E\rest{p^*}$. As above, for each $i<\theta$ let $g_i$ uniformise the restriction of~$\bar c$ to~$p_i$. We define $g\colon \alpha^*\to 2$  uniformising $\bar c$, as follows: for $\beta \in [\gamma_i,\gamma_{i+1})$, we let $g(\beta)= g_{i+1}(\beta)$, except that if $\beta\in E_{\alpha^*}$ we set $g(\beta) = c_{\alpha^*}(\beta)$. That~$g$ uniformises $\bar c$ follows from the fact that the changes from~$g_{i+1}$ are bounded: let $\beta\in p^*$, $\beta<\alpha^*$; so $\beta \in (\gamma_i,\gamma_{i+1})$ for some~$i$. Since $\cf(\beta)= \theta$, $E_{\alpha^*}\cap \beta$ is bounded below~$\beta$, so $g\rest{E_\beta} =^* g_{i+1}\rest{E_\beta} =^* f_\beta$. 

\smallskip

Suppose that the pair $(\task,\psi)$ is satisfied, say by $C$ and~$G$; then $W = W[G]\subset \Cof_\lambda(\theta)$ is stationary. The ladder system $\bar E\rest{W}$ has uniformisation: let $\bar c$ be a 2-colouring of $\bar E\rest{W}$;  let $H = G_{\bar c}\subseteq \SS^\task(C,G,\bar c)$ be $(C,G,\bar c)$-generic. Since each $p\in \QQ^\task(C)$ has uniformisation, $g_H = \bigcup H$ is a function with domain~$\lambda$ (for all $\gamma<\lambda$, the conditions in $\SS^\task$ with domain $\ge \gamma$ are dense), and it uniformises $\bar c$. 

\smallskip

It remains to show that $(\task,\psi)$ is nice. The argument follows \cite{Sh:64,Sh:186}. 

\medskip

After passing, possibly, to a generic extension, suppose that $S\subseteq \lambda$ is fat, and that  $\bar \PP = \seq{\PP_\zeta}_{\zeta\le \xi}$ is $S$-sparse and $\task$-friendly. We need to check that~$W = W[G_0]$ is stationary in $V(\PP_\xi)$. Let $D\in V(\PP_\xi)$ be a club. By extending by one step, we may assume that $\cf(\xi)\le \lambda$. 

\smallskip

Let~$\bar \tau$ be an associated sequence for $\bar \PP$. 
Let $Q$ be the collection of $\zeta<\xi$ for which $\QQ_\zeta = \SS^\task(S,G_0,\bar c)$ for some appropriate $\bar c \in V(\PP_\zeta)$; we write $\bar c^\zeta$ for~$\bar c$. By assumption, for $\zeta\in \xi\setminus Q$, if $\cf(\zeta)\le \lambda$ then $\s_\zeta$ is $W$-sparse; if $\cf(\zeta)> \lambda$ then we know that $\QQ_\zeta$ is $\lllambda$-closed. We further modify the construction of ~$\bar \tau$ to ensure that for all $\varsigma \le \xi$ with $\cf(\varsigma)\le \lambda$ (so $\tau_\varsigma$ is defined), for all $p\in \dom \tau_\varsigma$ and all $\zeta \in Q\cap \supp(p)$, $p\in A^\varsigma_\zeta$ and there is some string $\pi = \pi(p,\zeta)$ (in~$V$) such that $p\rest{\zeta}$ forces that $p(\zeta) = \pi$, and the length of~$\pi$ is $|\tau_\varsigma(p)|-1$. 

\smallskip

Elaborating only a little on our standard construction, we obtain a filtration $\bar N = \seq{N_i}_{i\le \theta}$ such that letting $\gamma_i = N_i\cap \lambda$, we have:
\begin{itemize}
	\item All objects above are elements of $N_0$, including an initial condition $p_{-1}\in \PP_\xi$;
	\item $\theta\subset N_0$, and if $\lambda = \mu^+$ then $\mu\subseteq N_0$;
	\item For limit $i\le \theta$, $\gamma_i \in S$;
	\item For successor $i<\theta$, $N_i^{<\theta}\subseteq N_i$.
\end{itemize}
We can do this since either~$\lambda$ is inaccessible, or $\lambda = \mu^+$ and $\mu^{<\theta} = \mu$. Note that the sequence $\seq{\gamma_i}$ is continuous. Let $\alpha^* = \gamma_\theta$. What is pertinent is that for successor $i<\theta$, $\cf(\gamma_i)\ge \theta$, so $E_{\alpha^*}\cap \gamma_i$ is bounded below~$\gamma_i$. Note that for limit $i<\theta$, $N_i^{<\theta}\subset N_{i+1}$. 

We use the trees-of-conditions method to obtain a sparse sequence $\bar p$ and an upper bound~$q$ forcing $\alpha^*\in W\cap D$; to do this, we need to take care of all possible choices for $p_i(\zeta)\rest{E_{\alpha^*}}$ for $\zeta\in Q\cap \supp q$. 

For $i\le \theta$ we define:
\begin{itemize}
	\item an ordinal $\delta_i$;
	\item sets $u_i\subset Q$; 
	\item functions $m_i \colon u_i \to i$.
\end{itemize}
 Let $T_i$ be the collection of all sequences $\bar f = \seq{f_\zeta}_{\zeta\in v}$, where $v\subseteq u_i$ is an initial segment of~$u_i$ and for all $\zeta\in v$, $f_\zeta \colon E_{\alpha^*}\cap \gamma_i\to 2$. We also define:
\begin{itemize}[resume]	
	\item for each $\bar f\in T_i$, a condition $p_i(\bar f)\in \PP_\xi$.
\end{itemize}
For a proper initial segment~$v$ of~$u_i$, let $\varsigma(v) = \min (u_i\setminus v)$; let $\varsigma(u_i)= \xi$. For $\bar f\in T_i$ defined on~$v$, we write $v(\bar f)=v$ and $\varsigma(\bar f) = \varsigma(v(\bar f))$. We say that $\bar f$ is \emph{maximal} (for~$T_i$) if $\dom \bar f = u_i$, i.e., if $\varsigma(\bar f) = \xi$. If $\bar f\in T_i$ and $\zeta \le \xi$ then we write $\bar f\rest{\zeta}$ for $\bar f\rest{(v\cap \zeta)}$. Note that whether maximal or not, $\cf(\varsigma(\bar f))\le \lambda$ (no $\zeta\in Q$ has cofinality $\lambda^+$), so in any case, $\tau_{\varsigma(\bar f)}$ is defined; we write $\tau_{v} = \tau_{\varsigma(v)}$ and $\tau_{\bar f} = \tau_{\varsigma(\bar f)} = \tau_{v(\bar f)}$. 

We ensure that the objects defined have the following properties:
\begin{equivalent}
	\item $u_i$, $m_i$, $T_i$ and the map $\bar f\mapsto p_i(\bar f)$ are all in $N_{i+1}$ (and in fact for successor~$i$, they will be in $N_i$);

	\item $|u_i|<\lambda$; if $\lambda = \mu^+$ then $|u_i|< \mu$;

	\item $u_{i}\subseteq u_j$ if $i<j$, and $u_j = \bigcup_{i<j} u_i$ for limit~$j$;

	\item if $i<j$ then $m_i = m_j\rest{u_i}$;

	\item $\delta_i = \gamma_i$ for limit~$i$; for successor~$i$,  $\gamma_{i-1}< \delta_{i} < \gamma_{i}$ and $\sup (E_{\alpha^*}\cap \gamma_i) < \delta_i$;

	\item For all $\bar f\in T_i$, $p_i(\bar f)\in \PP_{\varsigma(\bar f)}$, in fact $p_i(\bar f)\in \dom \tau_{\bar f}$, and $|\tau_{\bar f}(p_i(\bar f)| = \delta_i+1$;

	\item \label{item:conditions_trees:a_tree}
	For all $\bar f\in T_i$ and all $\zeta\in v(\bar f)$,  $p_i(\bar f)\rest{\PP_\zeta} = p_i(\bar f\rest \zeta)$; 

	\item \label{item:conditions_trees:something_else}
	For $\bar f\in T_i$, $v(\bar f)\subset \supp(p_i(\bar f))$; we write $\pi_i(\bar f,\zeta)$ for the string $\pi(p_i(\bar f),\zeta)$ mentioned above.

	\item \label{item:conditions_trees:guessing_f}
	For all $\bar f\in T_i$ and $\zeta\in v(\bar f)$, for all $\beta \in E_{\alpha^*} \cap [\gamma_{m_i(\zeta)}, \gamma_i)$, we have $\pi_i(\bar f,\zeta)(\beta) = f_\zeta(\beta)$. 

	\item \label{item:conditions_trees:bookkeeping_u}
	 For all $\bar f\in T_i$, $Q\cap \supp(p_i(\bar f))\subseteq u_\theta$;

	\item For all maximal $\bar f\in T_{i+1}$, $p_{i+1}(\bar f)$ forces that $D\cap [\gamma_i,\gamma_{i+1})$ is nonempty; 

	\item For all $j<i$, for all $\bar f\in T_i$, $p_i(\bar f)$ extends $p_j(\bar f[j])$, where $\bar f[j]\in T_j$ is the sequence $\seq{f_\zeta\rest{(E^*_\alpha\cap \gamma_j})}_{\zeta\in u_j \cap v(\bar f)}$;

	\item \label{item:conditions_trees:sparse}
	For all $j<i$ and $\bar f\in T_i$, if  $\varsigma(\bar f) = \varsigma(\bar f[j])$, then the sequence $\smallseq{p_k(\bar f[k]}_{k\in [j,i)}$ is sparse for $\tau_{\bar f}$;

	\item $0\in \supp (p_i(\bar f))$ for all $\bar f\in T_i$. 
\end{equivalent}

Let us show how to construct such objects. We start with $u_0 = \emptyset$; for $\bar f$ being the empty sequence (the only element of $T_0$) we let $p_0(\bar f)$ be some extension of $p_{-1}$ (the initial condition we started with) in $\dom \tau_\xi$; note that $\varsigma(\bar f) = \xi$. We let $\delta_0 = |\tau_\xi(p_0(\bar f))|-1$. We can ensure that $0\in \supp(p_0(\bar f))$. 

\smallskip

Suppose that $i\le \theta$ is a limit, and that all objects indexed by $j<i$ have been defined, and satisfy the properties above, except of course for \ref{item:conditions_trees:bookkeeping_u}. As required, we define $u_i = \bigcup_{j<i} u_j$, $m_i = \bigcup_{j<i} m_j$, and $\delta_i = \gamma_i$. Let $\bar f\in T_i$ and let $v = v(\bar f)$. There is some $j<i$ such that  $\varsigma(\bar f) = \varsigma(\bar f[j])$ (either $\varsigma(\bar f)\in u_j$, or it is~$\xi$); by \ref{item:conditions_trees:sparse}, we can let $p_i(\bar f) = \cb_{\tau_{\bar f}}(\smallseq{p_k[\bar f[k]]}_{k\in [j,i)})$. To define the sequence of objects up to~$i$, we need the sequence $\seq{E_{\alpha^*}\cap \gamma_j}_{j<i}$; as each $E_{\alpha^*}\cap \gamma_i$ has size $<\theta$, and $i<\theta$, this sequence is in $N_{i+1}$. Now \ref{item:conditions_trees:a_tree} for~$i$ follows from the canonical choices of bounds being, well, canonical, and $\tau_\zeta \cohere \tau_{\bar f}\rest{A^{\varsigma(\bar f)}_\zeta}$. \ref{item:conditions_trees:guessing_f} follows from $\pi_i(\bar f,\zeta) = \bigcup_{k\in [j,i)}\pi_j(\bar f[k],\zeta)$ (as $p_i(\bar f)$ extends $p_j(\bar f[j])$), where $j<i$ is any such that $\zeta\in u_j$. 

\smallskip

Suppose that $i<\theta$ is a successor ordinal, and that all objects have been defined for $j\le i-1$. First, we define $u_i$. This is done to make progress towards \ref{item:conditions_trees:bookkeeping_u}:
\begin{itemize}
	\item If $\lambda$ is inaccessible, then we can let $u_i = \bigcup_{\bar f\in T_{i-1}} Q\cap \supp (p_{i-1}(\bar f))$. 
	\item Otherwise, $\lambda = \mu^+$. If~$\mu$ is a limit cardinal, let $\seq{\mu_j}_{j<\theta}$ be a sequence of cardinals increasing to~$\mu$. In this case $|T_{i-1}|<\mu$, so we just ensure that for all $\bar f\in T_{i-1}$, ``$\mu_i$-much'' of $Q\cap \supp(p_{i-1}(\bar f))$ is added to~$u_i$ (note that $|\supp(p_{i-1}(\bar f))|$ is likely~$\mu$). 
	\item If $\mu = \nu^+$ is a successor, then we will likely have $|T_{i-1}| = \mu$ (even though $|u_{i-1}|\le \nu$), as $\theta = \mu$ and ${(<\theta)}^\nu = \mu$. We then add ``$i$-much'' of $Q\cap \supp(p_{i-1}(\bar f))$ for ``$i$-many'' $\bar f\in T_{i-1}$. 
\end{itemize}
We define $m_i$ to extend $m_{i-1}$ by letting $m_i(\zeta) = i-1$ for all $\zeta\in u_i \setminus u_{i-1}$. Note that~$u_i$ and~$\gamma_i$ determine~$T_i$. 

\smallskip

To define $p_i(\bar f)$ for $\bar f\in T_i$ we perform a transfinite construction. We work in~$N_{i}$. By our analysis immediately above, there is a regular cardinal $\kappa<\lambda$ with $ \kappa\ge  |T_i|$. Let $\bar M$ be a filtration (with all relevant objects in~$M_0$); we obtain a sequence $\seq{\epsilon_{\ell}}_{\ell <\kappa}$ with the usual properties: its limit points are in $S$, $M_{\epsilon_\ell}\cap \lambda = \epsilon_\ell$, and for $\ell<\kappa$ we have $\seq{\epsilon_{\ell'}}_{\ell'< \ell}\in M_{\epsilon_\ell+1}$. Since~$\kappa$ may be larger than~$\theta$, we cannot require that the sequence be continuous. 

Let $\seq{\bar f_\ell}_{\ell<\kappa}$ be a list of all maximal $\bar f\in T_i$, where each such $\bar f$ appears unboundedly often. This list is in~$M_0$. We also ensure that $\sup(E_{\alpha^*}\cap \gamma_i) <\epsilon_0$ (recall that $E_{\alpha^*}\cap \gamma_i$ is bounded below~$\gamma_i$ for successor~$i$). 

For $\ell \le \kappa$ and $\bar f\in T_i$ we define conditions $r_{\ell}(\bar f)$, satisfying the following:
\begin{orderedlist}
	\item $r_0(\bar f) = p_{i-1}(\bar f[i-1])$;

	\item Either $r_\ell(\bar f) = r_{0}(\bar f)$, or $v(\bar f)\subset \supp(r_\ell(\bar f))$ and $r_\ell(\bar f)\in \dom \tau_{\bar f}$, in fact  $r_\ell(\bar f)\in A^{\varsigma(\bar f)}_\zeta$  for all $\zeta\in v(\bar f)$;
\end{orderedlist}	

If $r_\ell(\bar f)\ne r_0(\bar f)$ then we write $\eta_\ell(\bar f) = |\tau_{\bar f}(r_\ell(\bar f))|$, and for $\zeta\in v(\bar f)$, we let $\pi_\ell(\bar f,\zeta) = \pi(r_\ell(\bar f),\zeta)$;

\begin{orderedlist}[resume]
	\item \label{item:partial_trees:zero}
	If $\zeta\in v(\bar f)$ then $r_\ell(\bar f\rest{\zeta})$ extends $r_\ell(\bar f)\rest{\PP_\zeta}$; 

	\item For $\ell<\kappa$, the map $\bar f\mapsto r_\ell(\bar f)$ is in $M_{\epsilon_\ell+1}$; 

	\item For $\ell'<\ell$, $r_{\ell}(\bar f)$ extends $r_{\ell'}(\bar f)$; 
\end{orderedlist}

We say that $r_{\ell}(\bar f)$ is \emph{new} if $\ell>0$ and for all $\ell'<\ell$, $r_{\ell'}(\bar f)\ne r_\ell(\bar f)$;

\begin{orderedlist}[resume]
	\item If $r_\ell(\bar f)$ is new then  $\eta_\ell(\bar f) \ge \epsilon_{<\ell} := \sup_{\ell'<\ell}\epsilon_{\ell'}$;


	\item \label{item:partial_trees:new_and_cohere}
	If $r_\ell(\bar f)$ is new then for all $\zeta\in v(\bar f)$, $r_{\ell}(\bar f\rest{\zeta})$ is new, and equals $r_{\ell}(\bar f)\rest{\PP_\zeta}$; 

	\item If $r_\ell(\bar f)\ne r_{0}(\bar f)$ then for all $\zeta\in v(\bar f)$, $\pi_\ell(\bar f,\zeta)$ extends~$f_\zeta\rest{[\gamma_{i-1},\gamma_i)}$.

	\item For all limit $\ell < \kappa$, the subsequence 
	\[
		\seq{ 
		r_{\ell'}(\bar f) \,:\, r_{\ell'}(\bar f)\text{ is new }
		}
	\]
	is sparse for $\tau_{\bar f}$ (but note that it may not be cofinal in $\smallseq{r_{\ell'}(\bar f)}_{\ell'<\ell}$, in which case the latter sequence is eventually constant.)
	\end{orderedlist}

For $\ell=0$ we follow \ref{item:partial_trees:zero}. At limit $\ell\le \kappa$ we let $r_\ell(\bar f)$ be either the eventually constant value of $r_{\ell'}(\bar f)$ for $\ell'<\ell$, if such exists; otherwise, we let $r_{\ell}(\bar f)$ be the $\tau_{\bar f}$-canonical upper bound of the sparse subsequence of new $r_{\ell'}(\bar f)$. Note that $r_\ell(\bar f)$ is new if and only if the second case holds, in which case, by \ref{item:partial_trees:new_and_cohere}, for all $\zeta\in v(\bar f)$, the second case holds for defining $r_{\ell}(\bar f\rest{\zeta})$, and it equals $r_\ell(\bar f)\rest{\PP_\zeta}$ by the coherence $\tau_\zeta \cohere \tau_{\bar f}\rest{A^{\varsigma(\bar f)}_\zeta}$. 

\smallskip

For the successor case, suppose that $r_{\ell}(\bar f)$ have been defined for all~$\bar f$. We now consider~$\bar f_\ell$ in steps:
\begin{enumerate}
	\item First, obtain the condition $s_0$ which is the ``sum'' of the sequence $\seq{r_{\ell}(\bar f_\ell\rest{v}}$ for~$v$ an initial segment of~$u_i$; for $\zeta\in \supp(s_0) = \bigcup_{v\le u_i}\supp(r_\ell(\bar f_\ell\rest{v}))$, having defined $s_0\rest{\zeta}$ extending each $r_\ell(\bar f_\ell\rest{v})\rest{\PP_\zeta}$, we define $s_0(\zeta)$ to be $r_\ell(\bar f_\ell\rest{\zeta+1})(\zeta)$. 

	Note that for all $\zeta \in u_i\cap \supp(s_0)$, $s_0\rest{\zeta}$ forces that $s_0(\zeta)$ is $\pi(s_0,\zeta)= \pi(r_\ell(\bar f_\ell\rest{\zeta+1}),\zeta)$, and either $|\pi(s_0,\zeta)|\le \gamma_{i-1}$, or $\pi(s_0,\zeta)$ extends $f_\zeta\rest{[\gamma_{i-1},\gamma_i)}$ (where $f_\zeta$ of course comes from the sequence $\bar f_\ell$). 

	\item Extend $s_0$ to a condition $s_1\in \PP_\xi$ by setting $\supp(s_1) = \supp(s_0)\cup u_i$; for each $\zeta\in u_i$ such that $|\pi(s_0,\zeta)|\le \gamma_{i-1}$, or $\zeta \notin \supp(s_0)$, we set $s_1(\zeta)$ to be some string (in~$V$) which extends $\pi(s_0,\zeta)$ if defined, and which extends $f_\zeta\rest{[\gamma_{i-1},\gamma_i)}$. 

	\item Extend~$s_1$ to a condition~$s_2$ in $\bigcap_{\zeta\in u_i}A^\xi_\zeta$, and also ensure that $s_2$ forces some $\beta\in [\gamma_{i-1},\gamma_i)$ into~$D$, where, recall, $D$ is the club in~$V(\PP_\xi)$ we want to get to meet~$W$ at~$\alpha^*$. Also ensure that $|\tau_\xi(s_2)|> \epsilon_{\ell-1}$. Further, by extending, we can ensure that for all initial segments~$v$ of $u_i$, $s_2\rest{\PP_{\varsigma(v)}}$ properly extends $s_1\rest{\PP_\zeta}$, and so properly extends $r_\ell(\bar f_\ell\rest{v})$. 

	\item Set $r_{\ell+1}(\bar f_\ell) = s_2$, and for all $\zeta\in u_i$, set $r_{\ell+1}(\bar f_\ell\rest{\zeta}) = s_2\rest{\PP_\zeta}$. For $\bar g\in T_i$ which is not an initial segment of $\bar f_\ell$, set $r_{\ell+1}(\bar g) = r_\ell(\bar g)$. 	
\end{enumerate}
Note that we ensured that $r_{\ell+1}(\bar g)$ is new if and only if $\bar g$ is an initial segment of $\bar f_\ell$. Also note that for $\zeta<\zeta'$ in~$u_i$, $r_{\ell+1}(\bar f\rest{\zeta'})\in A^{\zeta'}_\zeta$ follows from \cref{def:weak_coherent_sequence}\ref{item:weak_coherence:first_restriction}. 

\smallskip

This completes the construction of all $r_\ell(\bar f)$; we let $p_i(\bar f) = r_\kappa(\bar f)$. Since each $\bar f$ is tended to unboundedly many times, we see that each $r_\kappa(\bar f)$ is new; so for all~$\bar f$, $\eta_\kappa{\bar f} = \epsilon_{<\kappa}+1$, which we set to be $\delta_i+1$. This completes the construction of all $p_i(\bar f)$ for $i\le \theta$. 

\medskip

We want to find some maximal $\bar f\in T_\theta$ and some condition $q$, extending $p_i(\bar f)$ for all~$i<\theta$ (and as usual, not extending $p_\theta(\bar f)$), forcing that $\alpha^*\in W$; such a condition also forces that $\alpha^*\in D$. To this end, call a sequence $\bar f\in T_\theta$ \emph{acceptable} if there is some condition $q\in \PP_{\varsigma(\bar f)}$, extending $p_i(\bar f)$ for all $i<\theta$, which forces that $\alpha^*\in W$. So we want to show that some maximal $\bar f\in T_\theta$ is acceptable. 

To do that, we show by induction on initial segments $v$ of $u_\theta$, that some $\bar f$ with $v(\bar f) = v$ is acceptable. To do that, we show:
\begin{enumerate}
	\item The empty sequence is acceptable; and
	\item For all initial segments $w<u$ of~$u_\theta$, if $\bar g\in T_\theta$ with $v(\bar g)= w$ is acceptable, as witnessed by some $q$, then $\bar g$ can be extended to an acceptable~$\bar f$ with $v(\bar f) = u$, as witnessed by some $r$ extending~$q$. 
\end{enumerate}

Let us show~(1), that the empty sequence $\seq{}$ is acceptable. We have $0\in \supp(p_\theta(\seq{}))$, and $p_\theta(0)$ is a sequence of length $\alpha^*+1$ ending with~0. We let $q(0)$ agree with $p_\theta(\seq{})(0)$, except that we change the last~0 to a~1, i.e., we say $\alpha^*\in q(0)$. This condition extends $p_i(\seq{})$ for all $i<\theta$. Let $\upsilon = \varsigma(\seq{}) = \min u_\theta$; since $Q\cap \supp(p_\theta(\seq{}))\subseteq u_\theta$, we have $Q\cap \supp(p_\theta(\seq{})) = \emptyset$. Thus, we can define the condition~$q\in \PP_\upsilon$ by setting $\supp(q) = \supp(p_\theta(\seq{}))$, and by defining $q\rest{\zeta}$ by induction on $\zeta \le \upsilon$. This has already been done for $\zeta = 0$. If $\zeta>0$ is in $\supp(q)$ and $q\rest{\zeta}$ has already been defined, then as $\zeta\notin Q$, we know that $\QQ_\zeta$ is explicitly $W$-closed; $p_i(\seq{})(\zeta)$ is defined for a final segment of~$\zeta$, with height $\delta_i$; since $q\rest\zeta$ forces that $\alpha^*\in W$, it forces that there is some upper bound for $\seq{p_i(\seq{})}$, which we set to be $q(\zeta)$. 

\smallskip

Now we tend to~(2), which is proved by induction on (the order-type of) $u$. Suppose this has been proved for all initial segments $u'$ of~$u$. There are two cases, depending on the order-type of~$u$. First, suppose that~$u$ has a greatest element~$\upsilon$. Let $\varrho = \varsigma(u)$. By induction, it suffices to show~(2) for $w = u\cap \upsilon$. Suppose that $v(\bar g) = w$ and that~$\bar g$ is acceptable, as witnessed by some~$q$. 

We define $r\rest{\PP_\zeta}$ by induction on $\zeta\in [\upsilon,\varrho]$. We let $r\rest{\PP_\upsilon}$ be some extension of~$q$ which decides the value of $c^\upsilon_{\alpha^*}$ (where recall $\QQ_\upsilon = \SS^\task(S,G_0,\bar c^\upsilon)$), say it forces that~$c^\upsilon_{\alpha^*} = f_\upsilon$ (where $f_\upsilon\in V$); we let $\bar f = \bar g\conc f_\upsilon$. We also let $r(\upsilon) = \pi_\theta(\bar f,\upsilon)$; then $r\rest{\upsilon+1}$ is a condition since $\pi_\theta(\bar f,\upsilon)$ agrees with~$f_\upsilon$ from $\gamma_{m_\theta(\upsilon)}$ onwards. 

 We then repeat the argument for the empty sequence; we set $\supp(r)$ to agree with $\supp(p_\theta(\bar f))$ on $(\upsilon,\varrho)$, and note that there it is disjoint from~$Q$; for $\zeta\in (\upsilon,\varrho)\cap \supp(r)$ we let $r(\zeta)$ be an upper bound of $\seq{p_i(\bar f)(\zeta)}_{i\in [i^*,\theta)}$, which is forced by~$r\rest{\zeta}$ to exist since it forces that $\alpha^*\in W$. 

\smallskip
 
Finally, suppose that the order-type of~$u$ is a limit; let $\kappa<\lambda$  be regular and let $\seq{w_i}_{i\le \kappa}$ be an increasing and continuous sequence of initial segments of~$u$ with $w_\kappa = u$, with $w_0 = w$ being the initial segment we start with; let $\bar g$ be acceptable, as witnessed by some~$q_0$, with $v(\bar g) = w_0$. 

As expected, we work with a filtration $\bar N$ and a sequence $\seq{\epsilon_j}_{j<\kappa}$ with limit points in~$S$, $N_{\epsilon_j}\cap \lambda = \epsilon_j$, with $N_0$ containing all pertinent objects. We mimic the construction of \cref{lem:iterating_inverse_limit_case}.  We define sequences $\seq{q_j}_{j\le \kappa}$ of conditions and $\bar f_j$ satisfying:
\begin{sublemma}
	\item  $\seq{q_\ell, \bar f_\ell}_{\ell<j}, \in M_{\epsilon_j+1}$;
	\item $v(\bar f_j) = w_j$, and $q_j\in \PP_{\varsigma(w_j)}$ witnesses that $\bar f_j$ is acceptable; 
	\item If $\ell<j$ then $q_j$ extends~$q_\ell$ and $\bar f_\ell = \bar f_j\rest{w_\ell}$; 
	\item For successor $j<\kappa$, $q_j \in \dom \tau_{\varsigma(w_j)}$, indeed $q_j \in A^{\varsigma(w_j)}_\zeta$ for all $\zeta\in w_j$; and $|\tau_{\varsigma(w_j)}|> \epsilon_j$; 
	\item For all $\ell<j$, $\seq{q_k\rest{\PP_{\varsigma(w_\ell)}}}_{k\in (\ell,j)}$ is sparse for $\tau_{\varsigma{(w_\ell)}}$. 
\end{sublemma}
At successor steps we apply the induction hypothesis from $w_j$ to $w_{j+1}$, and then extend to a condition in $\dom \tau_{\varsigma(w_{j+1})}$ as required; at limit steps we take canonical sparse upper bounds and then an inverse limit. This completes the proof. 
\end{proof}


\bibliographystyle{halpha}
\bibliography{F1844_biblist}

\end{document}